\documentclass[a4paper,10pt]{amsart}

\usepackage{ifpdf}
\usepackage{amsmath, amsthm, amsfonts, amssymb, amscd}
\usepackage[active]{srcltx}
\usepackage{color}
\usepackage{mathtools}
\usepackage[ansinew]{inputenc}

\newtheorem{theorem}{Theorem}[section] %[section]
\newtheorem{corollary}[theorem]{Corollary}
\newtheorem{lemma}[theorem]{Lemma}

\newtheorem{proposition}[theorem]{Proposition}
\newtheorem*{proposition*}{Proposition}

\newtheorem*{question*}{Question}
\newtheorem*{theorem*}{Theorem}
\newtheorem*{claim*}{Claim}

\newtheorem*{corollary*}{Corollary}

\theoremstyle{definition}

\theoremstyle{remark}

\newtheorem*{remark*}{Remark}
\newtheorem*{definition*}{Definition}

\newtheorem*{remarks*}{Remarks}

\newcommand{\R}{\mathbb{R}}
\newcommand{\Z}{\mathbb{Z}}
\newcommand{\T}{\mathbb{T}}
\newcommand{\A}{\mathbb{A}}

\def\cotan{\qopname\relax o{cotan}}

\def\eqalignno#1{\displ@y \tabskip=\@centering
\halign to\displaywidth{\hfil$\@lign\displaystyle{##}$
\tabskip=0pt &$\@lign\displaystyle{{}##}$

\hfil\tabskip=\@centering
$\llap{$\@lign##$}\tabskip=Opt\crcr #1\crcr}}

\usepackage{graphicx}

\begin{document}

\author[P. Le Calvez]{Patrice Le Calvez}
\address{ Sorbonne Universit\'e, Universit\'e Paris-Cit\'e, CNRS, IMJ-PRG, F-75005, Paris, France \enskip \& \enskip Institut Universitaire de France}
\curraddr{}
\email{patrice.le-calvez@imj-prg.fr}

\title[Twist maps of the annulus: an abstract point of view]{Twist maps of the annulus: an abstract point of view}

\bigskip
\bigskip

\bigskip
\begin{abstract} We introduce the notion of abstract angle at a couple of points defined by two radial foliations of the closed annulus. We use this notion to give unified proofs of some classical results on area preserving positive twist maps of the annulus by using the Lifting Theorem and the Intermediate Value Theorem.  
 \end{abstract}
\maketitle



\maketitle

\bigskip
\bigskip

\bigskip
\bigskip

\section{Introduction}

\bigskip
Let us begin by introducing some notations and definitions. We will denote

\begin{itemize}
\item  $\T=\R/\Z$ the $1$-dimensional torus, 
\item $\A=\T\times[0,1]$ the closed annulus, 
\item  $\mathrm{int}(\A)=\T\times(0,1)$ the open annulus,
\item  $C_0=\T\times\{0\}$ and $C_1=\T\times \{1\}$ the boundary circles of $\A$,
\item  $\tilde \A=\R\times[0,1]$  the universal covering space  of $\A$,
\item   $\pi: (\tilde x,y)\mapsto (\tilde x+\Z,y)$  the covering projection,
\item $T:(\tilde x,y)\mapsto (\tilde x+1,y)$ the generating covering automorphism of $\pi$,
\item $\tilde X=\pi^{-1}(X)$ the preimage by the covering projection of a set  $X\subset \A$,
\item $\tilde d$ the Euclidean distance on $\tilde\A$,
\item $d$ the distance on $\A$ defined by $d(z,z')=\min_{\pi(\tilde z)=z, \,\pi(\tilde z')=z'}  \tilde d(\tilde z,\tilde z'\},$
\item $p_1: (x,y)\mapsto x$ and $p_2: (x,y)\mapsto  y$ the two projections defined on $\A$,
\item $\tilde p_1: (\tilde x,y)\mapsto \tilde x$ and $\tilde p_2: (\tilde x,y)\mapsto  y$ the two projections defined on $\tilde \A$,
\item $\omega=dx\wedge dy$ the usual area form on $\A$.
\end{itemize}

\medskip
Let $f$ be a homeomorphism of $\A$ isotopic to the identity and $\tilde f$ a lift of $f$ to $\tilde \A$ (such a lift commutes with $T$).

\begin{definition*} For every periodic orbit $O$ of $f$ of period $q\geq1$, there exists $p\in\Z$ such that for every $\tilde z\in \tilde O$ we have $\tilde f^q(\tilde z)=T^p(\tilde z)$. We will say that $O$ is a periodic orbit {\it of type $(p,q)$} for the lift $\tilde f$. 

\end{definition*}
\begin{definition*}A closed set $X\subset \A$ invariant by $f$ is {\it well ordered} if 
\begin{itemize}

\item 
$p_1{}_{\vert X}$ is injective,
 \item for every $\tilde z$, $\tilde z'$ in $\tilde X$ we have
$$\tilde p_1(\tilde z)<\tilde p_1(\tilde z')\Rightarrow \tilde p_1(\tilde f(\tilde z))<\tilde p_1(\tilde f(\tilde z')).$$
\end{itemize}
\end{definition*}
For example $C_0$ and $C_1$ are well-ordered invariant subsets. One deduces from the theory of homeomorphisms of the circle that if $X$ is a well ordered invariant closed set, there exists a real number $\rho_{\tilde f}(X)$, the {\it rotation number} of $X$, such that
for every
$\tilde z\in\tilde X$ and every $k\in\Z$, we have:
$$-1<\tilde p_1\circ \tilde f^k(\tilde z)-\tilde p_1(\tilde z)-k\rho_{\tilde f}(X)<1.$$ 
If $\rho_{\tilde f}(X)$ is rational and is written $\rho_{\tilde f}(X)=p/q$, with $p$ and $q$ relatively prime and $q>0$, then $X$ contains at least one periodic orbit of type $(p,q)$ for $\tilde f$.  If $\rho_{\tilde f}(X)$ is irrational, then $X$ contains a smallest closed invariant subset which is either an {\it invariant graph}, meaning the graph of a continuous function $\psi :\T^1\to[0,1]$, or a Cantor set.

\medskip
\begin{definition*}The map $f$ is a {\it positive twist map} if it is a $C^1$-diffeomorphism of $\A$ and if it satisfies the following equivalent properties:

\begin{enumerate}

\item  For every $\tilde x\in\R$, the map
$	y\mapsto \tilde p_{1}\circ \tilde f(\tilde x,y)$ is an increasing diffeomorphism from $[0,1]$ to its image  and $y\mapsto \tilde p_{1}\circ
\tilde f^{-1}(\tilde x,y)$ is an increasing diffeomorphism from $[0,1]$ to its image.

\item  The inequality $\displaystyle {\partial \,\tilde p_1\circ \tilde f\over \partial y} (\tilde x,y)>0$ holds for every $(\tilde x, y)\in \tilde\A$.

\item There exists $\beta\in]0,\pi/2[$ such that for every $z\in
\A$, the angle between the vertical vector $v=(0,1)$ and $Df(z).v$ belongs to $[-\beta, \beta-\pi]$
 and the angle between $v$ and $Df(z)^{-1}.v$ belongs to $[\pi-\beta, \beta]$.

\end{enumerate}
\end{definition*}

Twist maps, more particularly area preserving twist maps, have been widely studied in dynamical systems. They model the dynamics of a convex table billiard and appear naturally in the study of the dynamics of a generic area preserving diffeomorphism in a neighborhood of an elliptic fixed point. 
If $f$ is a positive twist map and $\tilde f$ is a lift of $f$, one proves easily that $\rho_{\tilde f}(C_0)\leq \rho_{\tilde f}(C_1)$, with a strict inequality if $f$ preserves the area. Let us state now three classical results about area preserving twist maps we will be interested in.

\begin{theorem}\label{th:birkhoffintro} Let $f$ be an area preserving positive twist map of $\A$ and $U$ an invariant set, homeomorphic to $\T\times[0,1)$ and containing $C_0$. Then there exists  a continuous map $\psi :\T\to(0,1]$ such that $U=\{(x,y)\in \A\,\vert \enskip 0\leq y<\psi(x)\}$. Moreover the natural lift $\tilde \psi:\R\to(0,1]$ of $\psi$ is $\cotan\beta$-Lipshitz.

\end{theorem}

\begin{theorem}  \label{th:aubrymatherintro}Let $f$ be an area-preserving positive twist map and $\tilde f$ a lift of $f$ to $\tilde \A$. For every $\rho\in[\rho_{\tilde f}(C_0),\rho_{\tilde f}(C_1)]$ there exists a well ordered invariant closed 
set
$X$ such that $\rho_{\tilde f}(X)=\rho$.
\end{theorem}

\begin{theorem}  \label{th:matherintro}Let $f$ be an area-preserving positive twist map of $\A$ with no invariant graph, except $C_0$ and $C_1$. Then there exist  $z$, $z'$ in $\mathrm{int}(\A)$ such that

$$\lim_{k\to-\infty} d(f^{k}(z), C_0)= \lim_{k\to+\infty} d(f^{k}(z), C_1)=0$$
and
$$ \lim_{k\to-\infty} d(f^{k}(z'), C_1)=\lim_{k\to+\infty} d(f^{k}(z'), C_0)=0.$$

\end{theorem}

Theorem \ref{th:birkhoffintro} is due to Birkhoff \cite{Bir1}, \cite{Bir2}. A rigorous exposition of Birkhoff arguments has been done by Herman  in \cite{He1}. This monography contains an appendix of Fathi \cite{Fa1} where an alternative proof is given using different topological arguments. One can also see Katznelson-Ornstein \cite{KatzO} or Siburg \cite{S}. Theorem \ref{th:aubrymatherintro} has been proved independently by Aubry \cite{AuL} and Mather \cite{Ma1}. The proofs use variational methods.
Nevertheless, one can obtain a proof using more topological arguments. One can prove that $p_1$ induces a $\cotan\beta$-bi-Lipshitz homeomorphism between a well ordered invariant closed set and its image by $p_1$. Using this fact Katok \cite{Kato} and Douady \cite{D} observed that to prove Theorem \ref{th:aubrymatherintro}, it is sufficient to show that for every rational number $p/q\in[\rho(C_0),\rho(C_1)]$, there exists a well ordered periodic orbit of type $p(p,q)$.  Moreover Bernstein \cite {Be} and Hall \cite{Hal} proved that for a positive twist map, the existence of a periodic orbit of type $(p,q)$ implies the existence of a well ordered periodic orbit of type $(p,q)$ (see Chenciner \cite{Ch} for a nice exposition of all these works). So, Theorem \ref{th:aubrymatherintro} can be deduced  from the classical Poincar\'e-Birkhoff Theorem:

\begin{theorem}  \label{th:poincarebirkhoffintro}Let $f$ be a homeomorphism of $\A$ isotopic to the identity and $\tilde f$ a lift of $f$ to $\tilde \A$. We suppose that $f$ preserves the area measure naturally defined by $\omega$ and that $\rho_{\tilde f}(C_0)<\rho_{\tilde f}(C_1)$. Then, for every rational number $p/q\in[\rho_{\tilde f}(C_0),\rho_{\tilde f}(C_1)]$ there exists at least one periodic orbit of type $(p,q)$.\end{theorem}

Fortunately, there exist simple proofs of Poincar\'e-Birkhoff Theorem in the particular case of a positive twist map, see Casdagli \cite{Ca} or \cite{L1}. In fact Poincar\'e-Birkhoff Theorem asserts that there exist at least two periodic orbits of type $(p,q)$. Similarly, in the statement of Theorem \ref{th:aubrymatherintro} one can add that $X$ contains at least two periodic orbits of type $(p,q)$ if $\rho=p/q$ is rational.

\medskip It has been proved by Moser \cite{Mo} that every positive twist map is the time-one map associated to a certain positive definite Lagrangian. This remark suggested Mather \cite{Ma3} that periodic positive definite Lagrangian systems should provide the right setting to extend the theory of positive twist maps to higher dimension, in particular in extending Theorem  \ref{th:aubrymatherintro} in a suitable way. Such systems have been widely studied, leading to the development of weak KAM theory (see Fathi \cite{Fa2}). Note that suitable extensions of Theorem \ref{th:birkhoffintro} have been stated by Arnaud \cite{Ar}, Bialy-Polterovich \cite{BiaP} and Herman \cite{He2}.

\medskip
 Coming back to the dimension two, it must be noticed that unlike variational methods, topological methods usually do not use the area-preserving condition but the fact $f$ {\it is not wandering}: for every non empty open set $U$, there exist $k>0$ such that $f^{-k}(U)\cap U\not=\emptyset $. In particular the four theorems above can be extended to this situation. Furthermore, 
the theorems above are still valid when $f$ is a finite composition of positive twist maps or more generally when $f$ is a{ \it positive tilt map}. It means that for every $x\in\T$ there exists a continuous negative function $\tau_x$ on $[0,1]$, satisfying $\tau_x(0)>-\pi/2$, such that for every $y\in[0,1]$, the angle between $v$ and $Df(x,y).v$ is equal to $\tau_x(y)$. Note that, in contrast with the set of positive twist maps, the set of tilt maps is stable by composition. In fact, informations about the angle of the tangents of the curves $f({\{x\}\times [0,1]})$ give informations about the vertical flips and the linking numbers of two distincts points (see Florio \cite{Fl} for related problems).

 \medskip
 If $f$ is a positive twist map, the negativity of the functions $\tau_x$  defined for a power $f^q$, $q\geq 1$, give informations on the shape of the curves $f^q(\{x\}\times[0,1])$. These informations, together with some results of plane topology,  are used in most of the topological proofs of the previous theorems, particularly in the seminal articles \cite{Bir2} and  \cite{Bir3}. Let $\mathcal V$ be the vertical foliation whose leaves are the segments $\{x\}\times [0,1]$ oriented upwards. One can find a negative function $\tau$ on $\A$, satisfying $\tau(z)>-\pi/2$ if $z\in C_0$, such that $\tau(z)$ measures the angle between the foliation $f^{-q}(\mathcal V)$ and the foliations $\mathcal V$ at $z$. This function can be lifted to a $T$-invariant function on $\tilde \A$, also denoted $\tau$. Its value at $\tilde z$ measures the angle between the lifted foliations $\tilde f^{-q}(\tilde{\mathcal V})$ and $\tilde{\mathcal V}$ at $\tilde z$, where $\tilde f$ is a lift of $f$. In the article, we will define an {\it abstract angle} in $\Z$, denoted $\tau(\tilde z, \tilde z', \mathcal F, \mathcal F')$ between two foliations at a couple of distincts points of $\tilde \A$. If $f$ is a positive tilt map and $q\geq 1$, the angle between $\tilde f^{-q}(\tilde{\mathcal V})$ and $\tilde{\mathcal V}$ at a couple $(\tilde z, \tilde z')$ will be non positive, and negative in case $\tilde z$ and $\tilde z'$ are on the same vertical of the lifted foliation. More precisely, if we set $\tilde W=\{(\tilde z, \tilde z')\in \tilde \A\times\tilde \A\,\vert\enskip \tilde z\not=\tilde z'\},$ we can define the angle $\dot\theta(\tilde z, \tilde z', \mathcal F)\in\Z/4\Z$ of a radial foliation $\mathcal F$ at $(\tilde z,\tilde z')\in \tilde W$. The set $\Z/4\Z$ is defined as the quotient space $(\R/2\pi\Z, \sim)$, where the equivalence classes of  $\sim$ are $\{0\}, (0,\pi), \pi, (\pi, 2\pi)$. It inherits a non Hausdorff topology but is path connected and its universal covering space is $\Z$, furnished with the  {\it digital line topology} also called {\it Khalimsky topology}. Denote $\mathfrak F$ the set of radial foliations. It inherits a natural topology and the map $\dot\theta: \tilde W\times \mathfrak F\to\Z/4\Z$ is continuous. Furthermore, $\mathfrak F$ is simply connected. Remind that a topological space $E$ is simply connected if it is path connected and if its fundamental group is trivial. The crucial point, immediate consequence of the Lifting Theorem, is that if $\tilde X\subset \tilde W$ is simply connected, the function $\dot\theta_{\vert\tilde X\times\mathfrak F}$ admits lift $\theta: \tilde X\times\mathcal F\to\Z$ and $$\tau(\tilde z, \tilde z', \mathcal F, \mathcal F')= \theta(\tilde z, \tilde z', \mathcal F')-\theta(\tilde z, \tilde z', \mathcal F).$$For this trivial geometry, every homeomorphism isotopic to the identity  is conformal. In particular if $\tilde X$ is invariant by $\tilde f\times \tilde f$ it holds that
 $\tau(\tilde f(\tilde z), \tilde f(\tilde z'), \mathcal V)\leq\theta(\tilde z, \tilde z', \mathcal V)$ with a strict inequality if $\tilde z$ and $\tilde z'$ are on the same vertical. We have got a Lyapounov function of $(\tilde f\times \tilde f)_{\vert \tilde X}$. We will see that natural invariant simply connected sets can be defined in each of our problems. The existence of a Lyapounov function will be extremely helpful. The pain is that we need to prove several formal results related to this abstract notion but the gain is that we have an efficient mathematical object to study positive twist maps. 
 
 \medskip 
We will give the detailed construction of abstract angles in the next section and will state several related properties. It must be highlighted that the construction works for homeomorphisms and so, we will introduce the notion of $\mathcal F$-monotone homeomorphism which generalizes the notion of tilt map. We will prove Theorem \ref{th:birkhoffintro}, Theorem \ref{th:matherintro} and the version of Theorem \ref{th:poincarebirkhoffintro} for twist maps in sections 3, 4, 5 respectively. More precisely we will state a version for $\mathcal F$-monotone homeomorphisms. The proofs are very close to the  classical ones but expressed in this new framework they show a lot of similarities by the use of the Lifting Theorem and the Intermediate Value Theorem.  In fact we will add a second proof of Theorem \ref{th:poincarebirkhoffintro}. It is very short and seems new, to the best of our knowledge.

\medskip
The new notion introduced in the article has been used by the author in a recent work about pseudo rotations \cite{L4}. There is no doubt that it can be used in some other situations concerning annular twist maps (for example Birkhoff attractors \cite{L3} or asymptotic linking numbers \cite {FlL}.  Can we go further? for example can it be used to obtain directly well ordered sets as stated in Theorem  \ref{th:aubrymatherintro}? can it be used to prove Theorem \ref{th:poincarebirkhoffintro} in its generality?

\section{Radial foliations and $\mathcal F$-monotone maps }\label{se:Calabi}  

\subsection{The space of radial foliation}
 
A {\it radial foliation} is an oriented topological foliation on $\A$ such that every leaf $\phi$ joins a point $\alpha(\phi)\in C_0$ to a point $\omega(\phi)\in C_1$ and satisfies $\phi\setminus\{\alpha(\phi), \omega(\phi)\}\subset \mathrm{int}(\A)$. We denote ${\mathfrak F}$ the set of radial foliations and $\mathrm{Homeo}_*(\A)$ the group of homeomorphisms of $\A$ furnished with the $C^0$ topology (which coincides with the uniform topology because $\A$ is compact). For every $\mathcal F\in\mathfrak F$ and every $f\in\mathrm{Homeo}_*(\A)$, we denote $f(\mathcal F)$ the radial foliation whose leaves are the images by $f$ of the leaves of $\mathcal F$. This action of $\mathrm{Homeo}_*(\A)$ on $\mathfrak F$ is transitive and the stabilizer $\mathrm{Homeo}_{\mathcal F}(\A)$ of $\mathcal F\in \mathfrak F$ is a closed subgroup of  $\mathrm{Homeo}_*(\A)$.  So, $\mathfrak F$ inherits a natural topology from the $C^0$ topology on $\mathrm{Homeo}_*(\A)$ and is simply connected for this topology. Indeed, $\mathfrak F$ is path connected and the morphism $i_* :\pi_1(\mathrm{Homeo}_{\mathcal F}(\A), \mathrm{Id})\to \pi_1(\mathrm{Homeo}_*(\A), \mathrm{Id})$  induced by the inclusion map $i:  \mathrm{Homeo}_{\mathcal F}(\A)\to \mathrm{Homeo}_*(\A)$ is bijective, both fundamental groups being infinite cyclic and generated by the same loop in $\mathrm{Homeo}_{\mathcal F}(\A)$ (see Hamstrom \cite{Ham}). In the whole article, the notation $\tilde{\mathcal F}$ means the lift of $\mathcal F\in\mathfrak F$ to the universal covering space $\tilde \A$. Note that every leaf $\tilde \phi$ of $\tilde{\mathcal F}$  is an oriented segment of $\tilde \A$ that separates $\tilde \A$ into two connected open sets, the component of  $\tilde \A\setminus\tilde\phi$ lying on the left of $\tilde \phi$ will be denoted $L(\tilde \phi)$ and the component lying on the right will be denoted $R(\tilde \phi)$. We get a total order $\preceq_{\mathcal F}$ on the set of leaves of $\tilde {\mathcal F}$ as follows:
$$\tilde \phi \preceq_{\mathcal F}  \tilde \phi'\Longleftrightarrow L(\tilde\phi)\subset L(\tilde\phi').$$
We can also define a partial order $\leq_{ {\mathcal F}}$ on $\tilde \A$:  write $\tilde z<_{{\mathcal F}}\tilde z'$ if $\tilde z$ and $\tilde z'$ are distinct and belong to the same leaf $\tilde\phi$ and if the segment of $\tilde\phi$ starting from $\tilde z$ and ending at $\tilde z'$ inherits the orientation of  $\tilde\phi$. For every $\tilde z\in \tilde\A$, we will denote $\tilde\phi_{\tilde z}$ the leaf of $\tilde{ \mathcal F}$ that contains $\tilde z$.

\subsection{Topological angles in the universal covering space}

Consider the natural projection
 $$\begin{aligned}\Pi: \Z&\to \Z/4\Z,\\
 k&\mapsto \dot k=k+4\Z.\end{aligned}$$ If $\Z/4\Z$ is endowed with the topology whose open sets are
 $$\emptyset, \{-\dot 1\}, \{\dot 1\}, \{-\dot 1, \dot 1\}, \{-\dot 1, \dot 0,\dot 1\},\{\dot 1, \dot 2, -\dot 1\}, \{-\dot 1,\dot 0, \dot 1, \, \dot 2\},$$ and $\Z$ with the topology generated
 by the sets $2k+1$ and $\{2k+1, 2k, 2k+1\}$, $k\in\Z$, then $\pi$ is a covering map. Note that both sets
$\Z/4\Z$ and $\Z$ are non Hausdorff but path connected. Note also that the connected subsets of $\Z$ are the intervals and consequently that every continuous function $\varphi: E\to \Z$ defined on a connected topological space $E$ satisfies the Intermediate Value Theorem.  Indeed, saying that $X\subset \Z$ is not an interval means that there exists $k\in\Z$ such that  $\left\{X\cap(-\infty,k], \,X\cap[k,\infty)\right\}$ is a partition of $X$. If $k$ is odd, it is a partition in open sets of $X$; if $k$ is even, it is a partition in closed sets of $X$. In both cases we deduce that $X$ is not connected.  Conversely suppose that $X$ is an interval and consider a non empty open and closed subset $Y$ of $X$. Let us prove that for every $k\in Y$, each point $k-1$ or $k+1$ belongs to $Y$ if it belongs to $X$. In the case where $k$ is even, this is true because $Y$ is open in $X$. In the case where  $k$ is odd, this is true because $X\setminus Y$ is open in $X$. One deduces that $Y=X$.  Note also that every continuous function $\varphi: E\to \Z$ defined on a compact topological space $E$ is bounded because the sets $\varphi^{-1} \left( (-2n,2n)\right)$, $n\geq 0$, are open. 

\medskip
Remind that we have defined
  $$\tilde W=\{(\tilde z, \tilde z')\in \tilde \A\times\tilde \A\,\vert\enskip \tilde z\not=\tilde z'\}.$$

We can define a continuous function  $\dot\theta :\tilde W\times{\mathfrak F}\to \Z/4\Z$  as follows:
$$\dot\theta(\tilde z, \tilde z', \mathcal F)= \begin{cases} -\dot  1 &\text{ if $\phi_{\tilde z} \prec_{\mathcal F} \phi_{\tilde z'}$,}\\
 \dot 0 &\text{ if $\tilde z<_{\mathcal F}\tilde z'$,}\\
\dot 1 &\text{ if $\phi_{\tilde z'} \prec_{\mathcal F} \phi_{\tilde z}$,}\\
 \dot 2 &\text{ if  $\tilde z'<_{\mathcal F}\tilde z$,}
\end{cases}$$
and we have
\begin{itemize}

\item $\dot\theta(\tilde z', \tilde z,\mathcal F) =\dot\theta(\tilde z, \tilde z',\mathcal F)+\dot 2$,

\item  $\dot\theta(\tilde f(\tilde z), \tilde f(\tilde z'),f(\mathcal F)) =\dot\theta(\tilde z, \tilde z',\mathcal F) $, for every  $f\in\mathrm{Homeo}_*(\A)$.
\end{itemize}

 For every $(\tilde z, \tilde z')\in \tilde W$, the function $\dot\theta_{\tilde z, \tilde z'}: \mathcal F\in{\mathfrak F}\to \dot\theta(\tilde z,\tilde z',\mathcal F)$ is continuous on $\mathfrak F$. This space being simply connected, the Lifting Theorem asserts that there exists a continuous function, 
$\theta: {\mathfrak F}\to \Z$, uniquely defined up to an additive constant in $4\Z$, such that $\Pi\circ \theta= \dot\theta_{\tilde z, \tilde z'}$. In particular, if $\mathcal F$ and $\mathcal F'$ are two radial foliations, the integer
$$ \tau(\tilde z, \tilde z', \mathcal F, \mathcal F') =\theta( \mathcal F')-  \theta( \mathcal F)$$
 does not depend on the choice of the lift $ \theta$.

Suppose that $f$ belongs to $\mathrm{Homeo}_*({\A})$ and that $\tilde f$ is a lift of $f$ to $\tilde \A$. The following results are immediate:
\begin{itemize}

\item  $\tau(\tilde z, \tilde z', \mathcal F, \mathcal F') = \tau(\tilde z', \tilde z, \mathcal F, \mathcal F')$,

\item $\tau(\tilde z,\tilde z', \mathcal F, \mathcal F')+\tau(\tilde z,\tilde z', \mathcal F', \mathcal F'')=\tau(\tilde z,\tilde z', \mathcal F, \mathcal F'')$,

\item $ \tau(\tilde z, \tilde z',  \mathcal F',  \mathcal F) =  -\tau(\tilde z, \tilde z', \mathcal F, \mathcal F')$,

\item  $\tau(\tilde f(\tilde z), \tilde f(\tilde z'),  f(\mathcal F),  f(\mathcal F')) =  \tau(\tilde z, \tilde z', \mathcal F, \mathcal F')$.

\end{itemize}

\subsection{Natural lifts }
By the Lifting Theorem, if $\tilde X$ is a simply connected subset of $\tilde W$, the map $\dot\theta_{\vert \tilde X\times \mathfrak F}$ admits a continuous lift $\theta:\tilde X\times \mathfrak F\to\Z $. Sometimes, a preferred lift is chosen, denoted $\theta^{\tilde X}$, that we call the {\it natural lift}. Of course, for every $(\tilde z, \tilde z')\in\tilde X$ the following holds:
$$ \tau(\tilde z, \tilde z', \mathcal F, \mathcal F') =\theta^{\tilde X}( \tilde z, \tilde z', \mathcal F')-  \theta^{\tilde X}( \tilde z, \tilde z', \mathcal F).$$Let us give some examples that will appear in the article.

\medskip
 If $\mathcal F\in\mathfrak F$, the set $\tilde X=\left\{(\tilde z,\tilde z')\in \tilde W\, \vert \enskip \dot\theta_{\tilde z,\tilde z'}(\mathcal F)\not=\dot 0\right\}$, is simply connected. The natural lift $\theta^{\tilde X}$ is the lift of $\dot\theta_{\vert \tilde X\times \mathfrak F}$ that is $\{1,2,3\}$-valued on $\tilde X\times\{\mathcal F\}$.

\medskip
 If $\mathcal F\in\mathfrak F$, the set $\tilde X=\left\{(\tilde z,\tilde z')\in \tilde W\, \vert \enskip \dot\theta_{\tilde z,\tilde z'}(\mathcal F)\in\{-\dot 1, \dot 0\}\right\}$, is simply connected. The natural lift $\theta^{\tilde X}$ is the lift of $\dot\theta_{\vert \tilde X\times \mathfrak F}$ that is $\{-1,0\}$-valued on $\tilde X\times\{\mathcal F\}$.

\medskip
The set $\tilde C_0\times\tilde C_1$ is simply connected. The natural lift $\theta^{\tilde C_0\times \tilde C_1}$  is the $\{-1,0,1\}$-valued lift of $\dot\theta_{\vert \tilde C_0\times\tilde C_1\times \mathfrak F}$. 

\medskip
More generally, if $\tilde X$ is a simply connected domain that contains $\tilde C_0\times \tilde C_1$, the  natural lift of  $\dot\theta_{\vert \tilde X\times \mathfrak F}$ is the $\Z$-valued lift that extends $\theta^{\tilde C_0\times \tilde C_1}$.

 \medskip
 
 Another case where a natural lift can be constructed is the following

\begin{itemize}
\item $\tilde X=\tilde X_0\times \tilde X_1$, where $\tilde X_0\cap \tilde X_1=\emptyset$,
\item every connected component of $\tilde X_0$ is simply connected and intersects $\tilde C_0$ in a non trivial connected set,
\item every connected component of $\tilde X_1$ is simply connected and intersects $\tilde C_1$ in a non trivial connected set.
\end{itemize}
The natural lift is the lift that is $\{-1,0,1\}$-valued on $(\tilde C_0\cap \tilde X_0)\times (\tilde C_1\cap \tilde X_1)\times \mathfrak F$.

\medskip
Let us present some other sets where a natural lift can be defined. We will denote $Y^c$ the complement of a subset $Y$ of a set $X$.

\medskip

Say that $U\subset\A$ is a {\it lower annulus} if it is open in $\A$, is homeomorphic to $\T\times[0,1)$ and contains $C_0$. The complement of $U$ is a connected closed set that contains
$C_1$.   Say that $U$ is {\it regular} if there exists an essential simple loop $\Gamma\subset \A$, meaning not freely homotopic to zero, such that $U$ is a connected component of $\Gamma^c$. Note that $\Gamma$ is the frontier of $U$. Upper annuli and regular upper annuli can be defined in the same way.

If $U$ is a lower annulus we set $\tilde U=\pi^{-1}(U)$. Note that if $U$ is regular, then $\tilde U\times \tilde U^c$ is simply connected and contains $\tilde C_0\times \tilde C_1$ and so the natural lift $\theta^{\tilde U\times \tilde U^c}$ is well defined.  Now, let us define $\theta^{\tilde U\times \tilde U^c}$ in the case where  $U$ is not necessarily regular. One can find a sequence $(U_n)_{n\geq 0}$ of regular lower annuli such that

\begin{itemize}

\item $\overline U_n\subset U_{n+1}$, for every $n\geq 0$,

\item   $U=\bigcup_{n\geq 0} U_n$.

\end{itemize} Observing that $\theta^{\tilde U_n\times \tilde U_n{}^c}$ and $\theta^{\tilde U_m\times \tilde U_m{}^c}$ coincide on $\tilde U_n\times \tilde U_m{}^c\times\mathfrak F$ if $n<m$, one deduces that  there exists a $\Z$-valued lift 
$\theta^{\tilde U\times \tilde U^c}$ of $\dot\theta_{\vert \tilde U\times \tilde U^c\times \mathfrak F}$, uniquely defined, that coincides with $\theta^{\tilde U_n\times \tilde U_n{}^c}$ on $\tilde U_n\times \tilde U^c\times \mathfrak F$, for every $n\geq 0$. Moreover $\theta^{\tilde U\times \tilde U^c}$ does not depend on the choice of the sequence $(U_n)_{n\geq 0}$. Note that the following properties  hold:
\begin{itemize} 
\item $\theta^{\tilde U\times \tilde U^c}$ is $\{-1,0,1\}$-valued on $\tilde U\times \tilde C_1\times \mathfrak F$;
\item if $U$ and $U'$ are lower annuli such that $U\subset U'$, then $\theta^{\tilde U\times \tilde U^c}$ and $\theta^{\tilde U'\times \tilde U'^c}$ coincide on $\tilde U\times \tilde U'{}^c\times\mathfrak F$;
\item if $U$ is a lower annulus, then for every $f\in \mathrm{Homeo}_*(\A)$, every lift $\tilde f$ of $f$ to $\tilde \A$, every $\tilde z\in \tilde U$, every $\tilde z'\in \tilde U^c$ and every $\mathcal F\in \mathfrak F$ it holds that 
$$\theta^{\tilde f(\tilde U)\times \tilde f(\tilde U)^c}(\tilde f(\tilde z), \tilde f(\tilde z'), f(\mathcal F))= \theta^{\tilde U\times \tilde U^c}(\tilde z, \tilde z', \mathcal F).$$
\end{itemize}
Of course, we have a similar version for upper annuli.

\medskip

 Say that an open set $U$ of $\A$ is a {\it lower disk} if it is open in $\A$, is homeomorphic to $\R\times[0,1)$ and meets $C_0$ (in that case $U\cap C_0$ is an open interval of $C_0$). The complement of $U$ in $\A$ is a connected closed set that meets $C_0$ and contains
$C_1$.  Say that $U$ is {\it regular} if one of the following situation occurs:
\begin{itemize}
\item $U$ is a connected component of $\gamma^c$, where $\gamma$ is a segment that joins a point of $ C_0$ to another point of $C_0$ and that does not meet $C_0$ elsewhere,
\item $U$ is a connected component of $\Gamma^c$, where $\Gamma$ is an essential loop that meets $C_0$ at a single point.
\end{itemize}
Upper disks and regular upper disks can be defined in the same way.

If $U$ is a regular lower disk, each connected component of $\tilde U=\pi^{-1}(U)$ is simply connected and intersects $\tilde C_0$ in a non trivial connected set. Moreover  $\pi^{-1}(U^c)=\tilde U^c$ is simply connected and contains $\tilde C_1$. So one can define a $\Z$-valued lift $\theta^{\tilde U\times \tilde U^c}$ of $\dot\theta_{\vert \tilde U\times \tilde U^c\times \mathfrak F}$ that coincides with $\theta^{\tilde C_0\times \tilde C_1}$ on $(\tilde U\cap \tilde C_0)\times (\tilde U^c\cap \tilde C_1)\times \mathfrak F$. By the same process as in the case of lower  annuli, this definition can be extended to general lower disks and we have the same properties as for lower annuli:
\begin{itemize} 
\item  if $U$ and $U'$ are lower disks such that $U\subset U'$, then $\theta^{\tilde U\times \tilde U^c}$ and $\theta^{\tilde U'\times \tilde U'^c}$ coincide on $\tilde U\times \tilde U'{}^c\times \mathfrak F$;

\item if $U$ is a lower disk, then for every $f\in \mathrm{Homeo}_*(\A)$, every lift $\tilde f$ of $f$ to $\tilde \A$, every $\tilde z\in \tilde U$, every $\tilde z'\in \tilde U^c$ and every $\mathcal F\in \mathfrak F$ it holds that 
$$\theta^{\tilde f(\tilde U)\times \tilde f(\tilde U)^c}(\tilde f(\tilde z), \tilde f(\tilde z'), f(\mathcal F))= \theta^{\tilde U\times \tilde U^c}(\tilde z, \tilde z', \mathcal F).$$
\end{itemize}
Here again, we have a similar version for upper disks.

\medskip

The following result, which admits a similar version for upper annuli and upper disks, will be useful in the article.
\begin{proposition} \label{pr:openannuli} Let $\mathcal F$ be a radial foliation of $\A$ and $U$ a lower annulus or a  lower disk. Define two subsets $K$ and $X$ of $U$ as follows:

\begin{itemize} 
\item $z\in K$ if there exist $ \tilde z\in \pi^{-1}(\{z\})$ and $\tilde z'\in \tilde U^c$ such that $ \theta^{\tilde U\times \tilde U^c}(\tilde z, \tilde z',\mathcal F)\geq 2$;
\item $z\in X$ if there exist $ \tilde z\in \pi^{-1}(\{z\})$ and $\tilde z'\in \tilde U^c$ such that $ \theta^{\tilde U\times \tilde U^c}(\tilde z, \tilde z',\mathcal F)= 2$.
\end{itemize}
Then, $K$ and $X$ coincide and  are closed in $U$.
\end{proposition}

\begin{proof}The inclusion $X\subset K$ is obviously true. Let us prove the converse.
We will apply the Intermediate Value Theorem. The problem is the fact that $\tilde U^c$ is not necessarily connected if $U$ is not regular and the Intermediate Value Theorem cannot be applied directly. Nevertheless we can avoid the difficulty by an approximation process. We consider a sequence $(U_n)_{n\geq 0}$ of regular lower annuli if $U$ is an annulus, of regular lower diks if $U$ is a disk, such $\overline U_n\subset U_{n+1}$, for every $n\geq 0$ and such that $U=\bigcup_{n\geq 0} U_n$. Fix $z\in K$. There exists $\tilde z\in \pi^{-1}(\{z\})$ and $ \tilde z'\in \tilde U^c $ such that $\theta^{\tilde U\times \tilde U^c}(\tilde z, \tilde z',\mathcal F)\geq 2$. There exists $n_0$ such that $\tilde z\in \tilde U_n$, for every $n\geq n_0$. Of course it holds that $ \tilde z'\in \tilde U_n^c $. The map $\theta^{\tilde U_n\times \tilde U_n{}^c}$ is continuous and $\{-1,0,1\}$-valued on $\tilde U_n\times \tilde C_1\times \mathfrak F$. By connectedness of $\tilde U_n^c$, one deduces that there exists $z''_n\in \tilde U_n ^c$ such that  $\theta^{\tilde U_n\times \tilde U_n{}^c}(\tilde z, \tilde z''_n,\mathcal F)=2$. The sequence $(\tilde z''_n)_{n\geq n_0}$ is bounded because for every $n\geq n_0$ it holds that  $\dot\theta_{\tilde z, \tilde z''_n}(\mathcal F)=\dot 2$. Consequently, taking a subsequence if necessary, one can always suppose that the sequence $(\tilde z''_n)_{n\geq 0}$ converges. If $\tilde z''$ is the limit, we have $\tilde z''\in \tilde U^c$ and $ \theta^{\tilde U\times \tilde U^c}(\tilde z, \tilde z'',\mathcal F)= 2$. We have proved that $z$ belongs to $X$.

\medskip
It remains to prove that $K=X$ is closed in $U$. It is not easy to prove directly that $K$ is closed and so we will consider the set $X$.
Let $(z_n)_{n\geq 0}$ be a sequence in $X$ that converges to $z\in U$. One can find a sequence $(\tilde z_n)_{n\geq 0}$ in $\tilde U$ and a sequence $(\tilde z'_n)_{n\geq 0}$ in $\tilde U^c$ such that for every $n\geq 0$, we have $$\pi(\tilde z_n)=z_n, \enskip  \theta^{\tilde U\times \tilde U^c}(\tilde z_n, \tilde z'_n,\mathcal F)=2.$$
Of course the choice of the lift $\tilde z_n$ is arbitrary and so, one can always suppose that the sequence  $(\tilde z_n)_{n\geq 0}$ is bounded. But this implies that the sequence $(\tilde z'_n)_{n\geq 0}$ is bounded because for every $n\geq 0$ it holds that  $\dot\theta_{\tilde z_n, \tilde z'_n}(\mathcal F)=\dot 2$. Consequently, taking a subsequence if necessary, one can always suppose that both sequences converge. The limit $\tilde z$ of the first one belongs to $\tilde U$ because it satisfies $\pi(\tilde z)=z$, the limit $\tilde z'$ of the second one belongs to $\tilde U^c$ because $\tilde U^c$ is closed. By continuity of $\theta^{\tilde U\times \tilde U^c}$ we deduce that $ \theta^{\tilde U\times \tilde U^c}(\tilde z, \tilde z',\mathcal F)= 2$. So, $z$ belongs to $X$.
\end{proof}

\subsection{ $\mathcal F$-monotone maps}

Let $\mathcal F$ be a radial foliation of $\A$. A map $f\in \mathrm{Homeo}_*(\A)$ is {\it  $\mathcal F$-increasing} if  for every $(\tilde z,\tilde z')\in \tilde W$ we have
$$\tau(\tilde z,\tilde z', \mathcal F, f^{-1}(\mathcal F))\geq 0$$ and
$$\tau(\tilde z,\tilde z', \mathcal F, f^{-1}(\mathcal F))=0\Longrightarrow \dot\theta_{\tilde z, \tilde z'}( \mathcal F)=\dot\theta_{\tilde z, \tilde z'}(f^{-1}( \mathcal F))\in\{-\dot 1, \dot 1\}.\footnote{ Note that the implication $\tau(\tilde z,\tilde z', \mathcal F, f^{-1}(\mathcal F))=0\Longrightarrow \dot\theta_{\tilde z, \tilde z'}( \mathcal F)=\dot\theta_{\tilde z, \tilde z'}(f^{-1}( \mathcal F))$ is obvious}$$
In particular it holds that
$$\dot\theta_{\tilde z, \tilde z'}(\mathcal F)\in\{\dot 0, \dot 2\}\Longrightarrow \tau(\tilde z,\tilde z', \mathcal F, f^{-1}(\mathcal F))>0\enskip \mathrm{and} \enskip \tau(\tilde z,\tilde z', f(\mathcal F), \mathcal F)>0.$$
The first inequality is obvious. To prove the second inequality, choose a lift $\tilde f$ of $f$ to $\tilde \A$ and
note that $$\begin{aligned}\
\dot\theta_{\tilde z, \tilde z'}(\mathcal F)\in\{\dot 0, \dot 2\}&\Longrightarrow  \dot\theta_{\tilde f^{-1}(\tilde z), \tilde f^{-1}(\tilde z')}(f^{-1}(\mathcal F))\in\{\dot 0, \dot 2\}\\
&\Longrightarrow  \tau(\tilde f^{-1}(\tilde z), \tilde f^{-1}(\tilde z'), \mathcal F, f^{-1}(\mathcal F))>0\\
&\Longrightarrow  \tau(\tilde z,\tilde z', f(\mathcal F), \mathcal F)>0.
\end{aligned}$$

Similarly, we will say that $f\in \mathrm{Homeo}_*(\A)$ is {\it  $\mathcal F$-decreasing} if  for every $(\tilde z,\tilde z')\in \tilde W$ we have
$$\tau(\tilde z,\tilde z', \mathcal F, f^{-1}(\mathcal F))\leq 0$$ and
$$\tau(\tilde z,\tilde z', \mathcal F, f^{-1}(\mathcal F))=0\Longrightarrow \dot\theta_{\tilde z, \tilde z'}( \mathcal F)=\dot\theta_{\tilde z, \tilde z'}(f^{-1}( \mathcal F))\in\{-\dot 1, \dot 1\}.$$

\medskip

What follows is a purely formal result.

\begin{proposition} \label{prop:formaltwist} Let $\mathcal F$ be radial foliation of $\A$, let $f$, $g$, $h$ be elements of $\mathrm{Homeo}_*(\A)$ and let $k$ be an integer.
\begin{enumerate}

\item The map $f$ is $\mathcal F$-increasing if and only if $f^{-1}$ is $\mathcal F$-decreasing.

\item If $f$ and $g$ are $\mathcal F$-increasing, then $f\circ g$ is $\mathcal F$-increasing.

\item If $f$ is $\mathcal F$-increasing, then $h\circ f\circ h^{-1}$ is $h(\mathcal F)$-increasing.

\item  If $f$ is $\mathcal F$-increasing, then $f$ is $f^k(\mathcal F)$-increasing.

\item If $f$ and $g$ are $\mathcal F$-decreasing, then $f\circ g$ is $\mathcal F$-decreasing.

\item If $f$ is $\mathcal F$-decreasing, then $h\circ f\circ h^{-1}$ is $h(\mathcal F)$-decreasing.

\item  If $f$ is $\mathcal F$-decreasing, then $f$ is $f^k(\mathcal F)$-decreasing.

\end{enumerate}
\end{proposition}

\begin{proof} Let us prove (1). It is sufficient to prove that $f^{-1}$ is $\mathcal F$-decreasing if $f$ is $\mathcal F$-increasing. Fix a lift $\tilde f$ of $f$ to $\tilde \A$. For every $(\tilde z,\tilde z')\in \tilde W$ it holds that
$$\begin{aligned} \tau(\tilde z,\tilde z', \mathcal F, f(\mathcal F))&= \tau(\tilde f^{-1}(\tilde z),\tilde f^{-1}(\tilde z'),f^{-1}(\mathcal F),  \mathcal F)\\
&= -\tau(\tilde f^{-1}(\tilde z),\tilde f^{-1}(\tilde z'), \mathcal F,  f^{-1}(\mathcal F))\\
&\leq 0.\end{aligned}$$

Moreover, we have
$$\begin{aligned} \tau(\tilde z,\tilde z', \mathcal F, f(\mathcal F))=0&\Longrightarrow \tau(\tilde f^{-1}(\tilde z),\tilde f^{-1}(\tilde z'), \mathcal F,  f^{-1}(\mathcal F))=0\\
&\Longrightarrow \dot\theta_{\tilde f^{-1}(\tilde z), \tilde f^{-1}(\tilde z')}( \mathcal F)=\dot\theta_{\tilde f^{-1}(\tilde z), \tilde f^{-1}(\tilde z')}(f^{-1}( \mathcal F))\in\{-\dot 1, \dot 1\}.\\
& \Longrightarrow\dot\theta_{\tilde z, \tilde z'}(f( \mathcal F))=\dot\theta_{\tilde z, \tilde z'}(\mathcal F)\in\{-\dot 1, \dot 1\}.
\end{aligned}$$

Now, let us prove (2). Suppose that $f$ and $g$ are $\mathcal F$-decreasing and fix a lift $\tilde g$ of $g$ to $\tilde \A$.  For every $(\tilde z,\tilde z')\in \tilde W$ it holds that
$$\begin{aligned} \tau(\tilde z,\tilde z', \mathcal F, g^{-1}\circ f^{-1}(\mathcal F))&= \tau(\tilde z,\tilde z', \mathcal F,  g^{-1}(\mathcal F))+\tau(\tilde z,\tilde z', g^{-1}(\mathcal F), g^{-1}\circ f^{-1}(\mathcal F))\\
&= \tau(\tilde z,\tilde z', \mathcal F,  g^{-1}(\mathcal F))+\tau(\tilde g(\tilde z), \tilde g(\tilde z'),\mathcal F, f^{-1}(\mathcal F))\\&\geq 0.
\end{aligned}$$

Moreover, if $\tau(\tilde z,\tilde z', \mathcal F, g^{-1}\circ f^{-1}(\mathcal F))=0$, then we have
$$\tau(\tilde z,\tilde z', \mathcal F,  g^{-1}(\mathcal F))=\tau(\tilde g(\tilde z), \tilde g(\tilde z'),\mathcal F, f^{-1}(\mathcal F))=0,$$ which implies that
$$\dot \theta_{\tilde z,\tilde z'} (\mathcal F) =   \dot\theta_{\tilde z, \tilde z'} (g^{-1}(\mathcal F))\in\{-\dot 1,\dot 1\}, \enskip \dot \theta_{\tilde g(\tilde z),\tilde g(z')} (\mathcal F) =   \dot\theta_{\tilde g(\tilde z), \tilde g(\tilde z')} (f^{-1}(\mathcal F))\in\{-\dot 1,\dot 1\}.$$

 Using the fact that $$\dot\theta_{\tilde z, \tilde z'} (g^{-1}(\mathcal F))=  \dot\theta_{\tilde g(\tilde z), \tilde g(\tilde z')} (\mathcal F), \enskip \dot\theta_{\tilde g(\tilde z), \tilde g(\tilde z')} (f^{-1}(\mathcal F))= \dot\theta_{\tilde z,\tilde z'} (g^{-1}\circ f^{-1}(\mathcal F)),$$ we deduce that
 $$\dot \theta_{\tilde z,\tilde z'} (\mathcal F) =  \dot \theta_{\tilde z,\tilde z'} (g^{-1}\circ f^{-1}(\mathcal F))\in\{-\dot 1,\dot 1\}.$$
 
 \medskip
To prove (3), suppose that  $f$ is $\mathcal F$-increasing. Fix a lift $\tilde h$ of $h$ to $\tilde \A$. For every $(\tilde z,\tilde z')\in \tilde W$ it holds that
$$\begin{aligned} \tau(\tilde z,\tilde z', h(\mathcal F), (h\circ f\circ h^{-1})^{-1}(h(\mathcal F)))&= \tau(\tilde z,\tilde z', h(\mathcal F),  h\circ f^{-1}(\mathcal F))\\
&= \tau(\tilde h^{-1}(\tilde z), \tilde h^{-1}(\tilde z'), \mathcal F,  f^{-1}(\mathcal F))\\&\geq 0.\end{aligned}$$

Moreover, if $\tau(\tilde z,\tilde z', h(\mathcal F), (h\circ f\circ h^{-1})^{-1}(h(\mathcal F)))=0$, then we have 
$$\tau(\tilde h^{-1}(\tilde z), \tilde h^{-1}(\tilde z'), \mathcal F,  f^{-1}(\mathcal F))=0,$$ and so we have
$$\dot \theta_{\tilde h^{-1}(\tilde z),h^{-1}(\tilde z')} (\mathcal F) =   \dot\theta_{\tilde h^{-1}(\tilde z), h^{-1}(\tilde z')} (f^{-1}(\mathcal F))\in\{-\dot 1,\dot 1\}.$$
It implies that
$$\dot \theta_{\tilde z,\tilde z'} (h(\mathcal F)) =   \dot\theta_{\tilde z, \tilde z'} (h\circ f^{-1}\ (\mathcal F))\in\{-\dot 1,\dot 1\}$$
which can be written
$$\dot \theta_{\tilde z,\tilde z'} (h(\mathcal F)) =   \dot\theta_{\tilde z, \tilde z'} ((h\circ f\circ h^{-1})^{-1}h(\mathcal F)))\in\{-\dot 1,\dot 1\}.$$

\medskip
To get (4), its is sufficient to apply (3) to $h=f^k$. The proofs of (5), (6), (7) are analogous to the proofs of (2), (3), (4) respectively.
\end{proof}

The following easy result will be useful in the article.

\begin{proposition} \label{prop:superiortwist}
Let $\mathcal F$ be a radial foliation of $\A$ and $f$ a $\mathcal F$-increasing homeomorphism of $\A$. If $(\tilde z_0,\tilde z_1)\in\widetilde W$ satisfies $\dot\theta_{\tilde z_0,\tilde z_1}(\mathcal F)= \dot 2$, then there exists a neighborhood $\tilde O_{\tilde z_0}$ of $\tilde z_0$ such that
$$\begin{aligned}
 \tilde z\in  \tilde O_{\tilde z_0} \enskip\mathit{and} \enskip \dot\theta_{\tilde z,\tilde z_1}(\mathcal F)= \dot 1\Longrightarrow \tau(\tilde z, \tilde z_1, \mathcal F, f^{-1}(\mathcal F))>1,\\ 
\tilde z\in  \tilde O_{\tilde z_0}\enskip\mathit{and} \enskip\dot\theta_{\tilde z,\tilde z_1}(\mathcal F)= -\dot 1\Longrightarrow \tau(\tilde z, \tilde z_1, \mathcal F, f(\mathcal F))<-1.\end{aligned}$$ 

\end{proposition}

\begin{proof} Consider the set $\tilde X=\left\{(\tilde z, \tilde z')\,\vert \enskip \dot\theta_{\tilde z, \tilde z'}(\mathcal F)\not=\dot 0\right\}$ and the natural lift $\theta^{\tilde X}$ of $\dot\theta_{\vert \tilde X\times \mathfrak F}$, which is $\{1,2,3\}$-valued on $\tilde X\times\{\mathcal F\}$. We have $\theta ^{\tilde X}(\tilde z_0, \tilde z_1,\mathcal F)=2$ and we deduce that  
 $$ \theta^{\tilde X}(\tilde z_0, \tilde z_1, f(\mathcal F))<2<\theta^{\tilde X}(\tilde z_0, \tilde z_1,f^{-1}(\mathcal F)).$$
So, there is a neighborhood $ \tilde O_{\tilde z_0}$ of $\tilde z_0$ such that for every $\tilde z\in  \tilde O_{\tilde z_0}$, we have $\tilde z\not=\tilde z_1$ and 
$$ \theta^{\tilde X}(\tilde z, \tilde z_1,f(\mathcal F))<2<\theta^{\tilde X}(\tilde z, \tilde z_1,f^{-1}(\mathcal F)).$$

If $\dot\theta_{\tilde z,\tilde z_1}(\mathcal F)= 1$, we have  $\theta^{\tilde X}(\tilde z, \tilde z_1,\mathcal F)=1$ and $\theta^{\tilde X}(\tilde z, \tilde z_1,f^{-1}(\mathcal F))>2$, which implies that $\tau(\tilde z, \tilde z_1, \mathcal F, f^{-1}(\mathcal F))>1$.

 If $\dot\theta_{\tilde z,\tilde z_1}(\mathcal F)= -1$, we have  $\theta^{\tilde X}(\tilde z, \tilde z_1,\mathcal F)=3$ and $\theta^{\tilde X}(\tilde z, \tilde z_1,f^{-1}(\mathcal F))<2$, which implies that $\tau(\tilde z, \tilde z_1, \mathcal F, f^{-1}(\mathcal F))<-1$.\end{proof}

We conclude this section with a result which permits to exhibit a natural class of $\mathcal F$-monotone maps. 

\begin{proposition} A positive twist map $f$ of $\A$ is $\mathcal V$-decreasing, where $\mathcal V$ is the vertical foliation whose leaves are the segments $\{x\}\times [0,1]$, $x\in\T$, oriented upwards. 

\end{proposition}

\begin{proof} By definition, if $f$ is a positive twist map, the following holds:
$$ \begin{aligned} \dot\theta_{\tilde z,\tilde z'} (\mathcal V)=\dot 0 &\Rightarrow  \dot\theta_{\tilde z,\tilde z'} (f^{-1}(\mathcal V))=-\dot 1,\\
 \dot\theta_{\tilde z,\tilde z'} (\mathcal V)=\dot 0 &\Rightarrow  \dot\theta_{\tilde z,\tilde z'} (f(\mathcal V))=\dot 1,\\
 \dot\theta_{\tilde z,\tilde z'} (\mathcal V)=\dot 2 &\Rightarrow  \dot\theta_{\tilde z,\tilde z'} (f^{-1}(\mathcal V))=\dot 1,\\
\dot\theta_{\tilde z,\tilde z'} (\mathcal V)=\dot 2 &\Rightarrow  \dot\theta_{\tilde z,\tilde z'} (f(\mathcal V))=-\dot 1.\end{aligned}$$
Applying the second implication to the couple $(\tilde f(\tilde z),\tilde f(\tilde z'))$, where $\tilde f$ is a lift of $f$, one obtains
$$ \dot\theta_{\tilde z,\tilde z'} (f^{-1}(\mathcal V))=\dot 0 \Rightarrow  \dot\theta_{\tilde z,\tilde z'} (\mathcal V)=\dot 1.$$
The set $\tilde X=\left\{(\tilde z, \tilde z')\,\vert \enskip \dot\theta_{\tilde z, \tilde z'}(\mathcal V)\in\{-\dot 1, \dot 0\}\right\}$ is simply connected and the natural lift $\theta^{\tilde X}$ of $\dot\theta_{\vert \tilde X\times \mathfrak F}$ is $\{-1,0\}$-valued on $\tilde X\times\{\mathcal V\}$.  The implication above tells us that $\dot\theta$ does not take the value $\dot 0$ on $\tilde X\times\{f^{-1}(\mathcal V)\}$, which implies that $\theta^{\tilde X}$ does not vanish on $\tilde X\times \{f^{-1}(\mathcal V)\}$.  

The set
 $\tilde Y=\left\{(\tilde z, T(\tilde z))\, \vert\enskip \tilde z\in \tilde \A\right\}$ is connected and included in $\tilde X$ and $\dot\theta$ is constant equal to $-\dot 1$ on $\tilde Y\times\mathfrak F$. The fonction $\theta^{\tilde X}$, being constant equal to $-1$ on $\tilde Y\times\{\mathcal V\}$, has the same property on $\tilde Y\times\mathfrak F$.  
The fonction $\theta^{\tilde X}$ does not vanish on $\tilde X\times \{f^{-1}(\mathcal V)\}$ and is constant equal to $-1$ on  $\tilde Y\times\{f^{-1}(\mathcal V)\}$. We deduce that it is negative on $\tilde X\times \{f^{-1}(\mathcal V)\}$. 

Consequently, the inequality $\tau(\tilde z, \tilde z', \mathcal V,  f^{-1}(\mathcal V))\leq 0$ holds if $\dot\theta_{\tilde z, \tilde z'}(\mathcal V)\in\{-\dot 1, \dot 0\}$ and is strict if $\dot\theta_{\tilde z, \tilde z'}(\mathcal V)=\dot 0$. In particular, it is also strict if $\dot\theta_{\tilde z, \tilde z'}(f^{-1}(\mathcal V))=\dot 0$. Using the fact that  $\tau(\tilde z', \tilde z, \mathcal V,  f^{-1}(\mathcal V))=\tau(\tilde z, \tilde z', \mathcal V,  f^{-1}(\mathcal V))$, we deduce that the inequality $\tau(\tilde z, \tilde z', \mathcal V,  f^{-1}(\mathcal V)\leq 0$ holds if $\dot\theta_{\tilde z, \tilde z'}(\mathcal V)\in\{\dot 1, \dot 2\}$ and is strict if $\dot\theta_{\tilde z, \tilde z'}(\mathcal V)=\dot 2$ or $\dot\theta_{\tilde z, \tilde z'}(f^{-1}(\mathcal V))=\dot 2$.\end{proof}

\section{Birkhoff Theory} \label{se:birkhofftheory}

Say that an essential simple loop $\Gamma\subset \A$ is {\it transverse} to a radial foliation $\mathcal F$ if it meets each leaf of $\mathcal F$ once. Similarly say that a simple path is transverse to $\mathcal F$ if it meets each leaf of $\mathcal F$ at most once.

Let $\varphi$ be a homeomorphism of a topological space $X$. Say that a subset $Y$ of $X$ is {\it wandering} if  the sets $\varphi^k(Y)$, $k\in\Z$, are pairwise disjoint. For instance, if $\varphi(Y)\subset Y$, then $Y\setminus \varphi(Y)$ is wandering. Remind that $\varphi$ is {\it non wandering} if there is no wandering open set but the empty set.

\medskip
The goal of this section is to prove the following result.
\begin{theorem} \label{th:birkhofftheorem} Let $\mathcal F$ be a radial foliation of $\A$ and $f$ a non wandering homeomorphism of $\A$ that is either $\mathcal F$-increasing or $\mathcal F$-decreasing. Then every invariant lower or upper annulus and every invariant lower or upper disk is regular and its frontier is transverse to $\mathcal F$.
\end{theorem}

\begin{proof} It is sufficient to prove the result in the case where $f$ is $\mathcal F$-increasing and $U$ is an invariant lower annulus or an invariant lower disk, the proof of the other cases being similar. Remind that $\theta^{ \tilde U\times \tilde U^c}$ is the  natural lift of $\dot\theta_{\vert \tilde U\times \tilde U^c\times \mathfrak F}$. Define three subsets $K$, $X$, $O$ of $U$ as follows:

\begin{itemize} 
\item $z\in K$ if there exist $ \tilde z\in \pi^{-1}(\{z\})$ and $\tilde z'\in \tilde U^c$ such that $ \theta^{\tilde U\times \tilde U^c}(\tilde z, \tilde z',\mathcal F)\geq 2$;
\item $z\in X$ if there exist $ \tilde z\in \pi^{-1}(\{z\})$ and $\tilde z'\in \tilde U^c$ such that $ \theta^{\tilde U\times \tilde U^c}(\tilde z, \tilde z',\mathcal F)=2$;
\item $z\in O$ if there exist $ \tilde z\in \pi^{-1}(\{z\})$ and $\tilde z'\in \tilde U^c$ such that $ \theta^{\tilde U\times \tilde U^c}(\tilde z, \tilde z',\mathcal F)>2$.

\end{itemize}

We have seen in Proposition \ref{pr:openannuli} that $K$ and $X$ coincide and are closed in $U$. Moreover the inclusion $O\subset K$ is obviously true.

\begin{lemma} \label{l:inclusionimage} We have 
$ f(X)\subset O$.
\end{lemma}
\begin{proof}For every $z\in f(X)$ there exist $z''\in \pi^{-1}(\{f^{-1}(z)\})$ and $\tilde z'\in \tilde U^c $ such that $ \theta^{\tilde U\times \tilde U^c}(\tilde z'', \tilde z',\mathcal F)=2$. We deduce that $\theta^{\tilde U\times \tilde U^c}(\tilde z'', \tilde z',f^{-1}(\mathcal F))>2$ because  $f$ is $\mathcal F$-increasing. The set $\tilde U$ being invariant by $\tilde f$, it means that $\theta^{\tilde U\times \tilde U^c}(\tilde f(\tilde z''), \tilde f(\tilde z'),\mathcal F)>2$, if $\tilde f$ is a given lift of $f$ to $\tilde \A$. Of course, it holds that $\pi(\tilde f(\tilde z''))=z$. Moreover we have $z\in U$ and $\tilde f(\tilde z')\in \tilde U^c$. So, $z$ belongs to $O$. 
\end{proof}

We will prove now that $K=\emptyset$ using the following sequence of assertions:

\begin{itemize}

\item $O$ is open because $\theta^{\tilde U\times \tilde U^c}$ is continuous and $\Z\cap(2,+\infty)$ is open;  

\item $O\setminus f(K)$ is wandering because
$f(K)=f(X)\subset O\subset K$; 

\item $f(K)$ is closed in $U$ because $f$ is a homeomorphism and $U$ is invariant;

\item $O\setminus f(K)$ is open because it is open in $U$ and $U$ is open;

\item $O=f(K)$ because  $O\setminus f(K)$ is open and wandering and $f$ is non wandering; 

\item $O$ is open and closed in $U$; 

\item $O=\emptyset$ because $O$ does not meet $C_0$ and 
$U$ is connected and meets $C_0$;

\item

$K=\emptyset$ because $f(K)\subset O$. 

\end{itemize}

We have proved that $ \theta^{\tilde U\times \tilde U^c}(\tilde z, \tilde z',\mathcal F)<2$ for every $\tilde z\in\tilde U$ and every $\tilde z'\in\tilde U^c$. We could prove similarly that $ \theta^{\tilde U\times \tilde U^c}(\tilde z, \tilde z',\mathcal F)>-2$. One deduces that
$$-2< \theta^{\tilde U\times \tilde U^c}(\tilde z, \tilde z',f^k(\mathcal F))<2$$ for  every $\tilde z\in\tilde U$, every $\tilde z'\in\tilde U^c$ and every $k\in\Z$, because  
$f$ is $f^k(\mathcal F)$-increasing. Consequently we have
$$ \begin{aligned}\theta^{\tilde U\times \tilde U^c}(\tilde z, \tilde z',\mathcal F)=1 &\Rightarrow  \theta^{\tilde U\times \tilde U^c}(\tilde z, \tilde z',f^{-1}(\mathcal F))=1,\\  \theta^{\tilde U\times \tilde U^c}(\tilde z, \tilde z',\mathcal F)=-1 &\Rightarrow  \theta^{\tilde U\times \tilde U^c}(\tilde z, \tilde z',f(\mathcal F))=-1,\end{aligned}$$
which implies that
$$ \begin{aligned}\dot\theta_{\tilde z, \tilde z'}(\mathcal F)=\dot 1 &\Rightarrow \tau(\tilde z, \tilde z', \mathcal F, f^{-1}(\mathcal F))=0,\\
\dot\theta_{\tilde z, \tilde z'}(\mathcal F)=-\dot 1 &\Rightarrow \tau(\tilde z, \tilde z', \mathcal F, f(\mathcal F))=0.\end{aligned}$$Note that for every leaf $\phi$ of $\mathcal F$, there exists a point $z_{\phi}\in U^c\cap \phi$, uniquely defined, such that for every $z\in \phi$, the conditions $z\in U$ and $\dot\theta_{z,z_{\phi}}(\mathcal F)=\dot 0$ are equivalent. By the Closed Graph Theorem, to conclude the proof of Theorem \ref{th:birkhofftheorem} it is sufficient to prove that $\{z_{\phi}\,\vert \enskip \phi\in \mathcal F\}$ is closed or equivalently that a point $z_0$ such that $\dot\theta_{z_0,z_{\phi}}(\mathcal F)=\dot 2$ does not belong to the closure of $U$. By Proposition \ref{prop:superiortwist}, if $\tilde z_{\phi}$ is a lift of $z_{\phi}$ and $\tilde z_0$ the lift of $z_0$ lying in the same leaf than $\tilde z_{\phi}$, there exists a neighborhood $\tilde O_{\tilde z_0}$ of $\tilde z_0$ such that
$$\begin{aligned}
 \tilde z\in \tilde O_{\tilde z_0}\enskip\mathrm{and} \enskip \dot\theta_{\tilde z,\tilde z_{\phi}}(\mathcal F)= \dot 1\Longrightarrow \tau(\tilde z, \tilde z_{\phi}, \mathcal F, f^{-1}(\mathcal F))>1,\\ 
\tilde z\in \tilde O_{\tilde z_0}\enskip\mathrm{and} \enskip\dot\theta_{\tilde z,\tilde z_{\phi}}(\mathcal F)= -\dot 1\Longrightarrow \tau(\tilde z, \tilde z_{\phi}, \mathcal F, f(\mathcal F))<-1.\end{aligned}$$ 
What has been proved above tells us that $\tilde O_{\tilde z_0}\cap\tilde U=\emptyset$. Indeed, a point $\tilde z\in \tilde O_{\tilde z_0}\cap\tilde U$ would satisfy $$\dot\theta_{\tilde z,\tilde z_{\phi}}( \mathcal F)=\dot 1\enskip \mathrm{and} \enskip \tau(\tilde z, \tilde z_{\phi}, \mathcal F, f^{-1}(\mathcal F))=0$$ or $$\dot\theta_{\tilde z,\tilde z_{\phi}}(\mathcal F)= -\dot 1\enskip\mathrm{and}\enskip \tau(\tilde z, \tilde z_{\phi}, \mathcal F, f(\mathcal F))=0.$$ \end{proof}

The first part of the proof of Theorem \ref{th:birkhofftheorem} can be deduced from the following result, which admits versions  for $\mathcal F$-decreasing homeomorphisms and upper annuli or upper disks. We will use this result in the next section. Remind that if $\varphi$ is a bijection of a set $X$, then $Y\subset X$ is {\it forward invariant } if $\varphi(Y)\subset Y$ and {\it backward invariant } if $Y\subset\varphi(Y)$.

\begin{proposition} \label{pr:halfbirkhofftheorem} Let $\mathcal F$ be a radial foliation of $\A$ and $f$ a $\mathcal F$-increasing and non wandering homeomorphism of $\A$. 
\begin{enumerate}
\item If $U$ is a backward invariant lower annulus or lower disk, we have\newline
$ \theta^{\tilde U\times \tilde U^c}(\tilde z, \tilde z',\mathcal F)<2$ for every $\tilde z\in \tilde U$ and $ \tilde z'\in \tilde U^c $.
\item If $U$ is a forward invariant lower annulus or lower disk, we have\newline
$ \theta^{\tilde U\times \tilde U^c}(\tilde z, \tilde z',\mathcal F)>-2$ for every $\tilde z\in \tilde U$ and $ \tilde z'\in \tilde U^c $.
\end{enumerate}
\end{proposition}

\begin{proof} We will prove (1), the proof of (2) being similar. Define the subsets $K$, $X$, $O$ of $U$ and the subset $O'$ of $f(U)$ as follows:

\begin{itemize} 
\item $z\in K$ if there exist $ \tilde z\in \pi^{-1}(\{z\})$ and $\tilde z'\in \tilde U^c$ such that $ \theta^{\tilde U\times \tilde U^c}(\tilde z, \tilde z',\mathcal F)\geq 2$;
\item $z\in X$ if there exist $ \tilde z\in \pi^{-1}(\{z\})$ and $\tilde z'\in \tilde U^c$ such that $  \theta^{\tilde U\times \tilde U^c}(\tilde z, \tilde z',\mathcal F)= 2$;
\item $z\in O$ if there exist $ \tilde z\in \pi^{-1}(\{z\})$ and $\tilde z'\in \tilde U^c$ such that $  \theta^{\tilde U\times \tilde U^c}(\tilde z, \tilde z',\mathcal F)>2$;
\item $z\in O'$ if there exist $ \tilde z\in \pi^{-1}(\{z\})$ and $\tilde z'\in \tilde f(\tilde U^c)$ such that \newline$  \theta^{\tilde f(\tilde U)\times \tilde f(\tilde U^c)}(\tilde z, \tilde z',\mathcal F)>2$.

\end{itemize}
We has seen that $K$ and $X$ coincide and are closed in $U$. Moreover we know that $O$ is open and included in $K$. We have the following generalization of Lemma  \ref{l:inclusionimage}.

\begin{lemma} \label{l:halfinclusionimage} We have $f(X)\subset O'$ and $ f(X)\cap U\subset O.$
\end{lemma}
\begin{proof} We need to slightly modify the proof of Lemma \ref{l:inclusionimage}. For every $z\in f(X)$ there exists $\tilde z''\in \pi^{-1}(\{f^{-1}(z)\})$ and $\tilde z'\in \tilde U^c $ such that $ \theta^{\tilde U\times \tilde U^c}(\tilde z'', \tilde z',\mathcal F)=2$. The map $f$ being $\mathcal F$-increasing we deduce that $ \theta^{\tilde U\times \tilde U^c}(\tilde z'', \tilde z',f^{-1}(\mathcal F))>2$.  This means that $ \theta^{\tilde f(\tilde U)\times \tilde f(\tilde U^c)}(\tilde f(\tilde z''), \tilde f(\tilde z'),\mathcal F)>2$, if $\tilde f$ is a given lift of $f$ to $\tilde \A$. We deduce that $z\in O'$.  Note that $\tilde f(\tilde z')\in \tilde U^c$ because $U^c$ is forward invariant. So, if  moreover we suppose that $z\in U$, then $\tilde f(\tilde z)$ belongs to $\tilde U$ and consequently 
$$ \theta^{ \tilde U\times  \tilde U^c}(\tilde f(\tilde z''), \tilde f(\tilde z'),\mathcal F)= \theta^{\tilde f(\tilde U)\times \tilde f(\tilde U^c)}(\tilde f(\tilde z''), \tilde f(\tilde z'),\mathcal F)>2.$$ So $z$ belongs to $O$. 
\end{proof}

Here again, we will prove that $K=\emptyset$ using the following sequence of assertions:
\begin{itemize}
\item $ (K\cup U^c)\setminus f(K\cup U^c)$ is wandering because 
$$f(K)\cap U\subset O\subset K, \enskip f(U^c)\subset U^c,$$
and so
$$\begin{aligned} f(K\cup U^c)&= f(K)\cup f(U^c)\\&\subset  f(K)\cup U^c \\
&=(f(K)\cap U)\cup U^c\\
&\subset O\cup U^c\\
&\subset K\cup U^c;
\end{aligned} $$

\item  $O\setminus (f(K)\cap U)$ is wandering because
 $$O\subset K, \enskip O\subset U, \enskip f(U^c)\subset U^c,$$ and so
 $$
 \begin{aligned}
 O\setminus (f(K)\cap U)&=(O\setminus f(K))\cup (O\setminus U)\\
 &=O\setminus f(K)\\
 &\subset (K\setminus f(K))\cap (K\setminus U^c)\\
 &\subset (K\setminus f(K))\cap (K\setminus f(U^c))\\
&=K\setminus f(K\cup U^c)\\
&\subset (K\cup U^c)\setminus f(K\cup U^c);\end{aligned}$$ 
\item $f(K)\cap U$ is closed in $U$ because $f(K)$ is closed in $f(U)$ and $U\subset f(U)$;

\item $O\setminus (f(K)\cap U)$ is open in $U$ because $O$ is open and $f(K)\cap U$ is closed in $U$;

\item $O\setminus (f(K)\cap U)$ is open because  $O\setminus (f(K)\cap U)$ is open in $U$ and $U$ is open;

\item $O=f(K)\cap U$  because $O\setminus (f(K)\cap U)$ is open and wandering and $f$ is non wandering;

\item  $O=\emptyset$ because $O$ does not meet $C_0$ and $U$ is connected and meets $C_0$; 

\item $O'=\emptyset$ because $f(U)$ is forward invariant and the previous arguments hold; 

\item $K=\emptyset$ because $f(K)=f(X)\subset O'$.

\end{itemize} 
\end{proof}
\begin{remarks*} 1. As explained above, there are analogous results for upper annuli or upper disks. For example, one can prove that if $U$ is a forward invariant upper annulus or upper disk, then
$\theta^{\tilde U^c\times \tilde U}(\tilde z, \tilde z',\mathcal F)>-2$ for every $\tilde z\in \tilde U^c$ and every $ \tilde z'\in \tilde U$.

\medskip
\enskip\enskip \enskip\enskip\enskip\enskip \, 2. Of course, Theorem \ref{th:birkhofftheorem} implies Theorem \ref{th:birkhoffintro}. If the hypothesis of Theorem \ref{th:birkhoffintro} are satisfied, one can apply Theorem \ref{th:birkhofftheorem} with the vertical foliation $\mathcal V$ and obtain a continuous function $\psi$ satisfying the first conclusion of Theorem \ref{th:birkhoffintro}. To prove that $\psi$ is Lipshitz it is sufficient to apply Theorem \ref{th:birkhofftheorem}  with the foliations $f^{-1}(\mathcal V)$ and $f(\mathcal V)$. It must be noticed that in all proofs of Theorem  \ref{th:birkhoffintro}, a forward (or backward) invariant set $K$ is defined and one proves that it is empty because the map is non wandering. The key point in the present formalism is that we have natural Lyapounov functions. More precisely,  if $\tilde X$ is a simply connected subset of $\tilde W$ invariant by $\tilde f\times \tilde f$ and if $\theta$ is a lift of $\dot\theta_{\vert \tilde X\times \mathfrak F}$, then the map $(\tilde z, \tilde z')\mapsto \theta(\tilde z, \tilde z', \mathcal F)$ does not decrease by the action of $\tilde f\times \tilde f$ and increases at $(\tilde z, \tilde z')$ if $\tilde z$ and $\tilde z'$ belong to the same leaf. We will repeat this argument in the rest of the article. 

\end{remarks*}

\section{Mather's theorem} \label{se:mathertheorem}

The aim of this section is to prove the following result, also true for $\mathcal F$-decreasing homeomorphisms.

\begin{theorem} \label{th:mathertheorem} Let $\mathcal F$ be a radial foliation of $\A$ and $f$ a $\mathcal F$-increasing and non wandering homeomorphism. If there is no invariant essential loop transverse to $\mathcal F$, but $C_0$ and $C_1$, then there exist $z$, $z'$ in $\mathrm{int}(\A)$ such that $$\lim_{k\to-\infty} d(f^{k}(z), C_0)= \lim_{k\to+\infty} d(f^{k}(z), C_1)=0$$
and
$$ \lim_{k\to-\infty} d(f^{k}(z'), C_1)=\lim_{k\to+\infty} d(f^{k}(z'), C_0)=0.$$

\end{theorem}

\begin{proof} We suppose that the hypothesis of Theorem  \ref{th:mathertheorem} are satisfied. We will prove that there exists a point $z$ such that $$\lim_{k\to-\infty} d(f^{k}(z), C_0)= \lim_{k\to+\infty} d(f^{k}(z), C_1)=0,$$
 the situation where $C_0$ and $C_1$ are exchanged being similar.

\begin{lemma} \label{le:continuabirkhofftheorem} If $X\subset\A$ is an invariant continuum, one of the following statements holds
$$X\cap (C_0\cup C_1)=\emptyset \,, \enskip C_0\cup C_1\subset X \,, \enskip X=C_0 \,, \enskip X=C_1.$$
\end{lemma}

\begin{proof} It is sufficient to prove that an invariant continuum $X$ that meets $C_1$ and does not contain $C_0$ is equal to $C_1$, because $C_0$ and $C_1$ have the same role. Consider $z\in X^c\cap C_0$. The set $X$ being a continuum, the connected component $U$ of $X^c$ that contains $z$ is either a lower annulus if $X\cap C_0=\emptyset$, or a lower disk if $X\cap C_0\not=\emptyset$. The map $f$ being non wandering and $X^c$ being open, there exists $p\geq 1$ such that $f^p(U)=U$ and such that the $f^k(U)$, $0\leq k<p$, are pairwise disjoint. Of course, we have $p=1$ if $U$ is an annulus. Let us prove by contradiction that $U$ is not a disk. If $U$ is disk, Theorem \ref{th:birkhofftheorem} tells us that the frontier of each $f^k(U)$, $0\leq k<p$, is an invariant transverse simple path. The union of these paths can be extended to an invariant transverse simple loop by adding $C_0\setminus \bigcup_{0\leq k<p}f^k(U)$. We have got a contradiction because this loop is different from $C_0$ and $C_1$. So $U$ is an annulus. Theorem \ref{th:birkhofftheorem} tells us that $U$ is regular and that its frontier is an invariant essential loop transverse to $\mathcal F$. By assumptions, it is equal to $C_1$, which means that $X=C_1$.
 
\end{proof}

\begin{corollary} \label{co:forwardcontinuabirkhofftheorem}  If $X_0\subset \A$ is a backward invariant continuum that strictly contains $C_0$ and does not meet $C_1$. Then we have 
$$\displaystyle \bigcap_{n\geq 0} f^{-n}(X_0)=C_0\,,\enskip \displaystyle C_1\subset \overline{\bigcup_{n\geq 0} f^n(X_0)}.$$ 
\end{corollary}

\begin{proof}Indeed, $ \bigcap_{n\geq 0} f^{-n}(X_0)$ and $\overline{\bigcup_{n\geq 0} f^n(X_0)}$ are invariant continua. The first one contains $C_0$ and does not meet $C_1$ while the second one strictly contains $C_0$. We conclude by Lemma \ref{le:continuabirkhofftheorem}.
\end{proof}

\begin{lemma} \label{le:continuabirkhofftheorembound} Suppose that:

\begin{itemize}
\item $X_0$ is a backward invariant continuum  such that $C_0\subset X_0$ and $X_0\cap C_1=\emptyset$;
\item$X_1$ is a forward invariant continuum such that $C_1\subset X_1$ and $X_1\cap C_0=\emptyset$;
\item we have $X_0\cap X_1=\emptyset$.
\end{itemize}
Then it holds that
\begin{itemize}
\item the connected component of $X_1{}^c$ that contains $X_0$, is a backward invariant lower annulus $U_0$; 
\item  the connected component of $X_0{}^c$ that contains $X_1$ is a forward invariant upper annulus $U_1$; 
 \item for every $\tilde z_0\in \tilde X_0$ and every $\tilde z_1\in \tilde X_1$ we have
 $$-2<\theta^{\tilde U_0\times \tilde U_0{}^c}(\tilde z_0, \tilde z_1,\mathcal F)=\theta^{\tilde U_1{}^c\times \tilde U_1}(\tilde z_0, \tilde z_1,\mathcal F)<2,$$ in particular it holds that $\dot\theta_{\tilde z_0,\tilde z_1}(\mathcal F)\not=\dot 2$.
\end{itemize}
\end{lemma}

\begin{proof} The two first items are obvious. To get the third one, we use Proposition \ref{pr:halfbirkhofftheorem} and the remark that follows this proposition. We know that
\begin{itemize}
\item 
$\theta^{\tilde U_0\times \tilde U_0{}^c}(\tilde z_0, \tilde z_1,\mathcal F)<2$ for every $\tilde z_0\in \tilde U_0$ and $ \tilde z_1\in \tilde U_0{}^c $.
\item 
$\theta^{\tilde U_1{}^c\times \tilde U_1}(\tilde z_0, \tilde z_1,\mathcal F)>-2$ for every $\tilde z_0\in \tilde U_1{}^c$ and $ \tilde z_1\in \tilde U_1 $.
\end{itemize}
It remains to prove that $\theta^{\tilde U_0\times \tilde U_0{}^c}$ and $\theta^{\tilde U_1{}^c\times \tilde U_1}$ coincide on $\tilde X_0\times \tilde X_1\times \{\mathcal F\}$.  They both lift $\dot\theta_{\vert \tilde X_0\times \tilde X_1\times \{\mathcal F\}}$ and coincide on $\tilde C_0\times \tilde C_1\times \{\mathcal F\}$. Moreover $\tilde X_0$ and $\tilde X_1$ are connected.  \end{proof}

\begin{lemma} \label{le:continuabirkhofftheoremintersection} Suppose that:

\begin{itemize}
\item $X_0$ is a backward invariant continuum  such that $C_0\varsubsetneq X_0$ and $X_0\cap C_1=\emptyset$,
\item$X_1$ is a forward invariant continuum such that $C_1\varsubsetneq X_1$ and $X_1\cap C_0=\emptyset$.
\end{itemize}
 Then there exists $n\geq 0$ such that $f^n(X_0)\cap X_1\not=\emptyset$. Moreover every point $z\in f^n(X_0)\cap X_1$ satisfies $$\lim_{k\to-\infty} d(f^{k}(z), C_0)= \lim_{k\to+\infty} d(f^{k}(z), C_1)=0.$$
 \end{lemma}

\begin{proof} We fix a lift $\tilde f$ of $f$ to $\tilde \A$. Consider $\tilde z_1\in \tilde X_1\setminus \tilde C_1$ and denote $\tilde z'_1$ the point of $\tilde C_1$ such that $\dot\theta_{\tilde z'_1, \tilde z_1}( \mathcal F)=\dot 2.$  By Proposition \ref{prop:superiortwist}, there exists a neighborhood $\tilde O_{\tilde z'_1}$ of $\tilde z'_1$ such that
$$\begin{aligned}
 \tilde z\in \tilde O_{\tilde z'_1}\enskip\mathrm{and} \enskip \dot\theta_{\tilde z,\tilde z_{1}}(\mathcal F)= \dot 1\Longrightarrow \tau(\tilde z, \tilde z_{1}, \mathcal F, f^{-1}(\mathcal F))>1,\\ 
\tilde z\in \tilde O_{\tilde z'_1}\enskip\mathrm{and} \enskip\dot\theta_{\tilde z,\tilde z_{1}}(\mathcal F)= -\dot 1\Longrightarrow \tau(\tilde z, \tilde z_{1}, \mathcal F, f(\mathcal F))<-1.\end{aligned}$$ 

 By Corollary \ref{co:forwardcontinuabirkhofftheorem}, if $n$ is large enough, one can find $\tilde z_0\in \tilde f^n(\tilde X_0)\cap \tilde O_{\tilde z'_1}$ such that $\dot\theta_{\tilde z_0,\tilde z_{1}}(\mathcal F)=\dot 1$. If $f^n(X_0)\cap X_1\not=\emptyset$, one can apply Proposition \ref{le:continuabirkhofftheorembound} with the foliations $\mathcal F$ and $f^{-1}(\mathcal F)$. We have 
 $$-2<\theta^{\tilde U_0\times \tilde U_0{}^c}(\tilde z_0, \tilde z_1,\mathcal F)=\theta^{\tilde U_1{}^c\times \tilde U_1}(\tilde z_0, \tilde z_1,\mathcal F)<2$$ and $$-2<\theta^{\tilde U_0\times \tilde U_0{}^c}(\tilde z_0, \tilde z_1,f^{-1}(\mathcal F))=\theta^{\tilde U_1{}^c\times \tilde U_1}(\tilde z_0, \tilde z_1,f^{-1}(\mathcal F))<2.$$
 We deduce that
 $$\theta^{\tilde U_0\times \tilde U_0{}^c}(\tilde z_0, \tilde z_1,f^{-1}(\mathcal F))=\theta^{\tilde U_0\times \tilde U_0{}^c}(\tilde z_0, \tilde z_1, \mathcal F)=1,$$ because $f$ is $\mathcal F$-increasing, which implies that
$$\tau(\tilde z_0, \tilde z_1, \mathcal F, f^{-1}(\mathcal F))=0.$$
We have found a contradiction.
\medskip
 
The fact that every point $z\in f^n(X_0)\cap X_1$ satisfies  $$\lim_{k\to-\infty} d(f^{k}(z), C_0)= \lim_{k\to+\infty} d(f^{k}(z), C_1)=0$$
 is an immediate consequence of  Corollary \ref{co:forwardcontinuabirkhofftheorem} that tells us that
$$\displaystyle \bigcap_{n\geq 0} f^{-n}(X_0)=C_0\,,\enskip \displaystyle \bigcap_{n\geq 0} f^{n}(X_1)=C_1.$$ 
\end{proof}

To conclude the proof of Theorem \ref{th:mathertheorem} it remains to prove that under the hypothesis of the theorem, there exists a backward invariant continuum $X_0$ that strictly contains $C_0$ and does not meet $C_1$. By similar arguments we could prove that there exists a forward invariant continuum that strictly contains $C_1$ and does not meet $C_0$. We will use an argument due to Birkhoff. Let $\Gamma$ be an essential simple loop that does not meet neither $C_0$ nor $C_1$. The component of $\Gamma^c$ that contains $C_0$  is a regular lower annulus $U$. We denote $X$ the connected component of $ \bigcap_{n\geq 0} f^{n}(\overline {U})$ that contains $C_0$. Of course $X$ is a backward invariant continuum. We will prove that $X\cap\Gamma\not=\emptyset$. Consider $n\geq 0$. One can find a regular lower annulus $U_n$ such that $f^k(U_n)\subset U$ if $0\leq k\leq n$. The open set $V_n=\bigcup_{n\geq 0} f^k(U_n)$ is backward invariant and the open set $W_n=\bigcup_{k\in\Z} f^k(U_n) =\bigcup_{k\leq 0} f^k(V_n)$ is invariant. The connected component of  $W_n{}^c$ that contains $C_1$ is an invariant continuum that does not contain $C_0$. So it is reduced to $C_1$. We deduce that $C_1\subset \overline {W_n}$. Note that $f(V_n)\setminus \overline {V_n}$ is an open wandering set, which implies that $f(V_n)\subset \overline V_n$. We deduce that $V_n$ and $W_n$ have the same closure. In particular, it holds that $C_1\subset\overline {V_n}$ and so, there exists $m>n$ such that $f^m(U^n)\not\subset U$ and $f^k(U^m)\subset U$ if $0\leq k<m$. One deduces that there exists a simple path $\gamma_n$ joining a point of $C_0$ to a point of $\Gamma$ that is contained in $f^m(U_n)\cap \overline {U}$. Note that $f^{-k}(\gamma_n)\subset \overline U$ if $0\leq k\leq n$. One can take a subsequence of the sequence $(\gamma_n)_{n\geq 0}$ that converges for the Hausdorff topology. Its limit is a continuum that meets $C_0$ and $\Gamma$ and that is contained in $\bigcup_{k\geq 0} f^k(\overline {U})$. We deduce that $X_0$ meets $\Gamma$. \end{proof}

\begin{remark*} Of course, Theorem \ref {th:mathertheorem} implies Theorem \ref {th:matherintro}. The proof given in this section is very close to the proof of  Theorem \ref {th:matherintro} given in \cite{L2}. In fact, Lemma \ref{le:continuabirkhofftheorem}, Corollary \ref{co:forwardcontinuabirkhofftheorem} and the final argument already appeared in \cite{L2}. Lemma
\ref{le:continuabirkhofftheoremintersection} also appeared but with a different proof, based on the semi-continuity of graphs naturally associated to the sets $X_0$ and $X_1$ and a continuation argument. The proof exposed in the present article is direct and simpler, based on Lemma  \ref{le:continuabirkhofftheorembound} which was not stated in \cite{L2}. 
\end{remark*}

\section{Poincar\'e-Birkoff Theorem} \label{se:poincarebirkhofftheorem}

%In this section, we denote $\theta^*:\tilde C_0\timdes \tilde C_1\tilde \mathfrak F\to\{-1,0,1\}$ the natural lift of $\dot\theta_{\tilde C_0\timdes \tilde C_1\tilde \mathfrak F}$.

In this section we will give two proofs of the following result.

\begin{theorem} \label{th:poincarebirkhofftheorem} Let $\mathcal F$ be a radial foliation of $\A$ and $f$ a $\mathcal F$-increasing and non wandering homeomorphism of $\A$. We suppose that a for certain lift $\tilde f$ of $f$ to $\tilde \A$, it holds that:

\begin{itemize}
\item for every $\tilde z\in  \tilde C_0$, we have $\tilde z\not =\tilde f(\tilde z)$ and  $\dot\theta_{\tilde z,\tilde f(\tilde z)}=-\dot 1$;
\item for every $\tilde z\in  \tilde C_1$, we have $\tilde z\not =\tilde f(\tilde z)$ and  $\dot\theta_{\tilde z,\tilde f(\tilde z)}=\dot 1$.
\end{itemize}

Then, there exist at least two  points of  $\A$ that are lifted to fixed points of $\tilde f$. 

\end{theorem}
\begin{proof}[First proof] 
We denote $X$ the set of points of $\A$ that are lifted to fixed points of $\tilde f$. We will argue by contradiction and suppose that $\# X\leq 1$.  We consider the topological sphere $S=\mathrm{int}(\A)\sqcup\left\{C_0, C_1\right\}$ obtained by crushing each boundary circle of $\A$ to a point. For every $\tilde z\in \pi^{-1}(\A\setminus X)$ and every $\mathcal F'\in\mathfrak F$ we have
$$\dot\theta _{T(\tilde z), \tilde f(T(\tilde z))} (\mathcal F')= \dot\theta _{\tilde z, \tilde f(\tilde z)} (\mathcal F').$$
So the map $(\tilde z, \mathcal F')\mapsto \dot\theta _{\tilde z, \tilde f(\tilde z)} (\mathcal F')$ lifts a continuous map $\delta :( \A\setminus X)\times\mathfrak F\to\Z/4\Z$, which itself lifts a continuous map $\dot\Delta : (S\setminus X)\times\mathfrak F\to\Z/4\Z$ defined as follows:
\begin{itemize}
\item $\dot \Delta(C_0, \mathcal F')= -\dot 1$,
\item $\dot \Delta(C_1, \mathcal F')= \dot 1$,
\item $\dot \Delta(z,\mathcal F')= \dot\theta _{\tilde z, \tilde f(\tilde z)} (\mathcal F')$ if  $z\in \mathrm{int}(\A)$ and $\pi(\tilde z)=z$.
\end {itemize}
By hypothesis on $X$, the set $S\setminus X$ is simply connected, and so $\dot\Delta$ admits a continuous $\Z$-valued lift $\Delta$. By connectedness of $\mathfrak F$, there exist $k_0\in\Z$ and $k_1\in \Z$ such that, for every $\mathcal F'\in\mathfrak F$, we have 
$$  \Delta(C_0, \mathcal F')= -1+4k_0\,, \enskip  \Delta(C_1, \mathcal F')= 1+4k_1.$$
The function $\Delta$ can be lifted to a continuous function $\delta: (\A\setminus X)\times\mathfrak F\to\Z$ equal to $-1+4k_0$ on $C_0\times\mathfrak F$ and to $1+4k_1$ on $C_1\times\mathfrak F$. We deduce from the following facts:
\begin{itemize}
\item$ \Pi(\delta(f(z), \mathcal F))=\dot\theta _{\tilde f(\tilde z), \tilde f^2(\tilde z)}(\mathcal F)= \dot\theta_{\tilde z, \tilde f(\tilde z)}(f^{-1}(\mathcal F))$ if $\pi(\tilde z)=z$,
\item$ \Pi(\delta(z, f^{-1}(\mathcal F))= \dot\theta_{\tilde z, \tilde f(\tilde z)}(f^{-1}(\mathcal F))$ if $\pi(\tilde z)=z$,
\item$ \delta(f(z), \mathcal F)= \delta(z, f^{-1}(\mathcal F))$ if $z\in C_0\cup C_1$,
\item $\A\setminus X$ is connected,
\end{itemize}
that the functions
$$z\mapsto \delta(f(z), \mathcal F)\,, \enskip z\mapsto \delta(z, f^{-1}(\mathcal F))$$ coincide on $\A\setminus X$.

Consequently,  the function $\delta_{\mathcal F}:z\mapsto \delta(z, \mathcal F)$ satisfies
$$\delta_{\mathcal F}(f(z))\geq \delta_{\mathcal F}(z)$$ and 
$$\delta_{\mathcal F}(f(z))>\delta_{\mathcal F} (z) \enskip\mathrm{if}\enskip\delta_{\mathcal F}(z)\in 2\Z.$$
By the Intermediate Value Theorem, every integer $k$ between $ -1+4k_0$ and $1+4k_1$ belongs to the image of $\delta_{\mathcal F}$. One can choose $k$ even. The set $K=\{z\in\tilde \A\,\vert \enskip \delta_{\mathcal F}(z)\geq k\}$ is closed and the set $O=\{z\in \tilde \A\,\vert \enskip \delta(z)>k\}$ is open. Moreover, we have $f(K)\subset O\subset K$. We deduce that $O\setminus f(K)$ is a non empty wandering open set. We have found a contradiction. 
 \end{proof}

\begin{proof}[Second proof]  We begin with a lemma:

\begin{lemma} \label{l:lecalvezlemma} Let $\mathcal F_0$ and $\mathcal F_1$ be two radial foliations of $\A$. Let $\tilde\phi_0$ be a leaf of $\tilde{\mathcal F_0}$  and $\tilde\phi_1$ be a leaf of $\tilde{\mathcal F_1}$. We suppose that $\tilde\phi_0\cap\tilde\phi_1\not=\emptyset$ and we define $\tilde z_0$ and $\tilde z_1$ in $\tilde \phi_0\cap \tilde\phi_1$ by the following properties
$$\begin{aligned}\tilde z\in (\tilde\phi_0\cap\tilde\phi_1)\setminus \{\tilde z_0\}\Longrightarrow 
\dot\theta_{\tilde z_0,\tilde z}(\mathcal F_0)=\dot 0\\
\tilde z\in (\tilde\phi_0 \cap\tilde\phi_1)\setminus \{\tilde z_1\}\Longrightarrow 
\dot\theta_{\tilde z_1,\tilde z}(\mathcal F_1)=\dot 2.
\end{aligned}$$
Then if $\tilde z_0\not=\tilde z_1$, we have $\tau(\tilde z_0, \tilde z_1, \mathcal F_0, \mathcal  F_1)=0$.
\end{lemma}

\begin{proof}  We suppose that the hypothesis of the lemma are satisfied and that $\tilde z_0\not=\tilde z_1$. Consider the segments
$$\tilde V_0=\{\tilde z\in \tilde \A\setminus \{\tilde z_0\} \,\vert\enskip \dot\theta_{\tilde z_0,\tilde z} (\mathcal F_0)=\dot 2\}\cup\{\tilde z_0\}, \enskip \tilde V_1=\{z\in \tilde \A\setminus \{\tilde z_1\} \,\vert\enskip \dot\theta_{\tilde z_1,\tilde z}( \mathcal F_1)=\dot 0\}\cup\{\tilde z_1\}.$$ 
Note that $\tilde X_0=\tilde C_0\cup\tilde V_0$  and $\tilde X_1=\tilde C_1\cup\tilde V_1$ are disjoint and that the function $\dot\theta$ does not take the value $\dot 2 $ on 
$\tilde X_0\times \tilde X_1\times\{\mathcal F_0,\mathcal F_1\}.$ The set
$\tilde X_0\times \tilde X_1$ being simply connected,  one can define the natural lift  $\theta^{\tilde X_0\times \tilde X_1}$ of $\dot\theta_{\vert \tilde X_0\times \tilde X_1\times\mathfrak F}$. It is $\{-1,0,1\}$-valued on $\tilde X_0\times \tilde X_1\times\{\mathcal F_0,\mathcal F_1\}.$ In particular we have
$$\tau(\tilde z_0, \tilde z_1, \mathcal F_0, \mathcal  F_1)=\theta^{\tilde X_0\times \tilde X_1}(\tilde z_0,\tilde z_1,\mathcal F_1)-\theta^{\tilde X_0\times  \tilde X_1}(\tilde z_0,\tilde z_1,\mathcal F_2)=0-0=0.$$
\end{proof}
 
Let us prove first that for every $\tilde\phi\in\tilde{\mathcal F}$ we have $\tilde f^{-1}(\tilde\phi)\cap \tilde\phi\not=\emptyset$. The map $\tilde z\in \tilde\phi \mapsto \dot\theta_{\tilde z,\tilde f(\tilde z)}(\mathcal F)$ takes the value $-\dot 1$ at the intersection point of $\tilde C_0$ and $\tilde\phi$ and takes the value $\dot 1$ at the intersection point of $\tilde C_1$ and $\tilde\phi$. By connectedness of $\tilde\phi$, it takes the value $\dot 0$ or $\dot 2$ at a point $\tilde z\in \tilde\phi$, which implies that $\tilde f(\tilde z)$ belongs to $\tilde\phi$. 

Moreover, the map $f$ being $\mathcal F$-increasing, if two different points $\tilde z$ and $\tilde z'$ satisfies
$\dot\theta_{\tilde z,\tilde z'}(\mathcal F)=\dot\theta_{\tilde z,\tilde z'}(f^{-1}(\mathcal F))\in\{\dot 0, \dot 2\}$, then we have  $\tau(\tilde z, \tilde z', \mathcal F, f^{-1}( \mathcal  F))\not=0$. Consequently, if $\tilde \phi$ is a leaf of $\tilde{\mathcal F}$, the two points $\tilde z_0$ and $\tilde z_1$ defined in the statement of Lemma \ref{l:lecalvezlemma} coincide if we set $\mathcal F_0=\mathcal F$, $\mathcal F_1=f^{-1}(\mathcal F)$, $\tilde \phi_0=\tilde\phi$ and $\tilde \phi_1=\tilde f^{-1}(\tilde\phi)$.  We set $\tilde z_{\tilde\phi}=\tilde z_0=\tilde z_1$. By definition, if $\tilde z\in \tilde\phi\cap \tilde f^{-1}(\tilde \phi)$ is different from $\tilde z_{\tilde\phi}$ we have 
$$\dot\theta_{\tilde z_{\phi},\tilde z}(\mathcal F)=\dot 0\, ,\enskip
\dot\theta_{\tilde z_{\phi},\tilde z}(f^{-1}(\mathcal F))=\dot 2,$$
or equivalently
$$\dot\theta_{\tilde z_{\phi},\tilde z}(\mathcal F)=\dot 0\, ,\enskip
\dot\theta_{\tilde f(\tilde z_{\phi}),\tilde f(\tilde z)}(\mathcal F)=\dot 2.$$
For analogous reasons, there exists a point $\tilde z'_{\tilde\phi}\in  \tilde\phi\cap \tilde f^{-1}(\tilde \phi) $ such that for every $z\in \tilde\phi\cap \tilde f^{-1}(\tilde \phi)$ different from  $\tilde z'_{\tilde\phi}$, we have 
$$\dot\theta_{\tilde z_{\phi},\tilde z}(\mathcal F)=\dot 2\, ,\enskip
\dot\theta_{\tilde z'_{\phi},\tilde z}(f^{-1}(\mathcal F))=\dot 0,$$
or equivalently
$$\dot\theta_{\tilde z'_{\phi},\tilde z}(\mathcal F)=\dot 2\, ,\enskip
\dot\theta_{\tilde f(\tilde z'_{\phi}),\tilde f(\tilde z)}(\mathcal F)=\dot 0.$$

The set 
$$\begin{aligned}X&=\left\{z\in \A\,\vert \enskip \pi(\tilde z)=z\enskip\Rightarrow \tilde\phi_{\tilde f(\tilde z)}=\tilde\phi_{\tilde z} \right\}\\&=\left\{z\in \A\,\vert \enskip \pi(\tilde z)=z\enskip\mathrm{and} \enskip f(\tilde z)\not=\tilde z\Rightarrow \enskip  \dot\theta_{\tilde z, \tilde f(\tilde z)}(\mathcal F)\in\{\dot 0,\dot 2\}  \right\}\end{aligned}$$ is closed and separates $C_0$ and $C_1$. So there exists at least one connected component of $X$ that separates $C_0$ and $C_1$ (for instance see Newman \cite{N}, Theorem 14.8, page 123). 
For every leaf $\tilde \phi$ of $\tilde{\mathcal F}$, the sets
$$\{z\in\A\,\vert \enskip \pi(\tilde z)=z\Rightarrow \enskip \tilde z\not=\tilde z_{\tilde\phi}\enskip \mathrm {and} \enskip \dot\theta_{\tilde z, \tilde z_{\tilde\phi}}(\mathcal F)=\dot 0\}$$ and 
$$\{z\in \A\,\vert \enskip \pi(\tilde z)=z\Rightarrow \enskip z\not=\tilde z_{\tilde\phi}\enskip \mathrm {and}\enskip \dot\theta_{\tilde z, \tilde z_{\tilde\phi}}(f^{-1}(\mathcal F))=\dot 2\}$$ are disjoint  from $X$. One deduces that $\pi(\tilde z_{\tilde\phi})$ belongs to every separating connected component of $X$. In particular there exists a unique separating component $Y$ of $X$. For the same reasons every $\pi(\tilde z'_{\phi})$ belongs to $Y$.
Define $$Y_0= \left\{z\in Y\,\vert \enskip \pi(\tilde z)=z\enskip\mathrm{and} \enskip f(\tilde z)\not=\tilde z\Rightarrow \enskip\dot\theta_{\tilde z, \tilde f(\tilde z)}(\mathcal F)=\dot 0\right\}$$ and
$$Y_2= \left\{z\in Y\,\vert \enskip \pi(\tilde z)=z\enskip\mathrm{and} \enskip f(\tilde z)\not=\tilde z\Rightarrow\enskip  \dot\theta_{\tilde z, \tilde f(\tilde z)}(\mathcal F)=\dot 2 \right\}.$$
Both sets are closed and their union separates $C_0$ and $C_1$. If their intersection is not connected, there are at least two points $z\in Y$ such that $\tilde f(\tilde z)=\tilde z$ if $\pi(\tilde z)=z$ and we are done. If their intersection is connected, it is known that at least one of the sets $Y_0$ or $Y_1$ separates $C_0$ and $C_1$ (for instance see Newman \cite{N}, Theorem 16.2, page 127).  Let us suppose for example that $Y_0$ satisfies this property. One can find a connected component $Y'_0$ of $Y_0$ that separates $C_0$ and $C_1$. As explained above, all points $\pi(\tilde z_{\phi})$ and $\pi(\tilde z'_{\phi})$, $\tilde \phi\in\tilde{\mathcal F}$, belong to $Y'_0$. One can suppose that there is at most one point in $Y'_0$ that is lifted to fixed points of $\tilde f$, otherwise we are done. We will suppose that such a point $z_*$ exists, the case where there is no such point being a little bit easier with a similar proof.  Note that for every leaf $\tilde\phi$ and every point $\tilde z\in \tilde Y'_0$ belonging to $\tilde\phi$ and different from $\tilde z'_{\tilde\phi}$ we have $$\dot\theta_{\tilde z, \tilde z'_{\tilde\phi}}(\mathcal F)=\dot 0\, ,\enskip
\dot\theta_{\tilde f(\tilde z), \tilde f(\tilde z'_{\tilde \phi})}(\mathcal F)=\dot 2.$$
It implies that every point $\tilde z\in \pi^{-1}(\{z_*\})$ coincide with $\tilde z'_{\tilde \phi}$ for a certain leaf $\tilde\phi$ and moreover that $$\pi(\tilde z'_{\tilde\phi})\not=z_*\Rightarrow\dot\theta_{\tilde z'_{\phi},\tilde f(\tilde z'_{\tilde\phi})}(\mathcal F)=\dot 0.$$ 
Consequently, the sets $Y'_0$ and $f(Y'_0)$ have a unique intersection point $z_*$ and the set $f(Y'_0)\setminus\{z_*\}$ is contained in the component of $Y'_0{}^c$ that contains $C_1$. This contradicts the fact that $f$ is non wandering.
\end{proof}

\begin{remark*} The second proof given in this section is very close to the proof of Poincar\'e-Birkhoff Theorem for positive twist maps given in \cite{L1}. What follows Lemma \ref{l:lecalvezlemma} was already done in \cite{L1} except that the existence of the second fixed point was not shown. Lemma \ref{l:lecalvezlemma} was stated only in the case where $\mathcal F_0=\mathcal V$ and $\mathcal F_1=f^{-q}(\mathcal V)$, the map $f$ being a positive twist map, and its proof was based on an induction process. The proof we present here is much simpler. The very short proof given before could have be written using Euclidean angles, instead of integer angles but with a clumsier redaction. Amazingly, it has never be done, as far as I know. 
\end{remark*}

%\ographystyle{gstart}s
%\renewcommand{\refname}{\centerline{\Large \bf  Bibliographie}}

%\Addresses

\end{document}